\documentclass[a4paper,oneside,11pt]{article}
\usepackage{amsmath}
\usepackage{amsfonts}
\usepackage{amssymb}
\usepackage{titlesec}
\usepackage{graphicx}
\usepackage{fancyhdr}
\usepackage{float}
\usepackage{pdfsync}
\usepackage{latexsym}
\usepackage{mathrsfs}
\usepackage{booktabs}
\usepackage{mathtools}
\usepackage{amsthm}
\usepackage{wasysym}
\usepackage{cite}
\usepackage{color}
\usepackage{geometry}
\usepackage{extarrows}
\usepackage{enumerate}
\usepackage{dsfont}
\usepackage{bm}
\usepackage[dvips]{epsfig}
\usepackage{amscd}
\usepackage[marginal]{footmisc}
\usepackage[bookmarks,colorlinks,citecolor=blue,linkcolor=red]{hyperref}
\allowdisplaybreaks

\def\QEDopen{{\setlength{\fboxsep}{0pt}\setlength{\fboxrule}{0.2pt}\fbox{\rule[0pt]{0pt}{1.3ex}\rule[0pt]{1.3ex}{0pt}}}} 
\def\QED{\QEDopen} 

\def\endproof{\hspace*{\fill}~\QED\par\endtrivlist\unskip}

\numberwithin{equation}{section}

\newtheorem{theorem}{Theorem}[section]
\newtheorem{lem}[theorem]{Lemma}
\newtheorem{prop}[theorem]{Proposition}
\newtheorem{rem}{Remark}[section]
\renewcommand{\thefootnote}{}

\title{ Large Time Behavior of Solutions to Cauchy Problem for 1-D Compressible Isentropic Navier-Stokes/Allen-Cahn System }
\date{  }
\author{Yazhou C{\small HEN}$^1$, Qiaolin H{\small E}$^2$, Xiaoding S{\small HI}$^{1*}$ \\[3mm]
\scriptsize$^{1}$ {College of Mathematics and Physics, Beijing University of
Chemical Technology, Beijing 100029, China}\\
\scriptsize$^{2}$ {School of Mathematics, Sichuan University, Chengdu, Sichuan, 610065,  China} }  

\begin{document}
\maketitle
\renewcommand{\thefootnote}{\fnsymbol{footnote}}
\footnotetext[1]{{Corresponding author. }\\
{Email: chenyz@mail.buct.edu.cn (Y.Chen), qlhejenny@scu.edu.cn (Q.He), shixd@mail.buct.edu.cn (X.Shi) }}
\begin{abstract}
This paper is concerned with the large time behavior of the solutions to the Cauchy problem for the one-dimensional compressible Navier-Stokes/Allen-Cahn system with the immiscible two-phase flow initially located near the phase separation state. Under the assumptions that the initial data is a small perturbation of the constant state, we prove the global existence and uniqueness of the solutions and establish the time decay rates of the solution as well as its higher-order spatial derivatives. Moreover, we derive that the solutions of the system are time asymptotically approximated by the solutions of the modified parabolic system and obtain decay rates in $L^2$ and $L^1$. Furthermore, we show that the solution of the system is time asymptotically approximated in $L^p (1 \leq p \leq+\infty)$ by the diffusion waves.

\end{abstract}
\noindent{\bf keywords:} Navier-Stokes/Allen-Cahn equation, well-posedness, large-time behavior, spectral analysis, time decay rates

\

\noindent{\bf AMS subject classifications:} 35Q35; 35B65; 76N10; 35M10; 35B40; 35C20; 76T30

\
\section{Introduction}
\indent\qquad
The coupled compressible Navier-Stokes/Allen-Cahn (called as NSAC) system  is an important diffusion interface model for describing compressible immiscible two-phase flow, where Navier-Stokes equation is used to simulate the motion of the mixture and Allen-Cahn equation is used to represent the changes of the phase field variable, which is widely used in chemical materials, petroleum exploitation and other fields.
The 1-D compressible barotropic NSAC system in the Eulerian coordinates  is as follows (see Blesgen \cite{b1999}, Heida-M\'{a}lek-Ralagopal \cite{HMR2012}):
\begin{equation}\label{nsac-e}
\left\{\begin{array}{llll}
\displaystyle  \rho_{\tilde{t}}+( \rho u)_{\tilde{x}}=0,\\
\displaystyle  \rho u_{\tilde{t}}+ \rho uu_{\tilde{x}}+P( \rho)_{\tilde{x}}=(\tilde{\nu}(\rho) u_{\tilde{x}})_{\tilde{x}}-\frac{\epsilon}{2}(\chi_{\tilde{x}}^2)_{\tilde{x}}, \\
\displaystyle  \rho\chi_{\tilde{t}}+ \rho u\chi_{\tilde{x}}=-\mu,\\
\displaystyle  \rho\mu= \frac{\rho}{\epsilon}(\chi^3-\chi)-\epsilon \chi_{\tilde{x}\tilde{x}},
\end{array}\right.
\end{equation}
for $(\tilde{x},\tilde{t})\in\mathds{R}\times\mathds{R}_+$.
The unknown functions $ \rho, u,$ and $\chi$ denote the total density, the mean velocity, and the concentration difference of the mixture fluids, respectively. $\mu$ is the chemical potential and $P=P( \rho)$ is the pressure, 
$\epsilon>0$ is the thickness of the interface between the phases, $\tilde{\nu}=\tilde{\nu}(\rho)>0$ the viscosity coefficient and  is assumed to be a smooth function of $\rho$.
The initial conditions are following
\begin{equation}\label{init-e}
 ( \rho,u,\chi)(\tilde{x},0)=(\rho_0,u_0,\chi_0)(\tilde{x}), 
\end{equation}
where $\rho_0>0$. As the space variable tends to infinity, we assume
\begin{equation}\label{inf-e}
\lim_{x\rightarrow\pm\infty}( \rho_0, u_0, \chi_0)(\tilde{x})= ( \bar{\rho}, \bar{u}, \bar{\chi}),
\end{equation}
where  $ \bar{\rho}>0, \bar{u}$  are the given constants and $\bar{\chi}=\pm 1$.
To be clear,  the function value of $\chi$ can be used to determine the area occupied by the two phases of the fluid,  where  $\{(\tilde{x},\tilde{t})\big|\chi(\tilde{x},\tilde{t}) = 1\}$ is occupied by fluid 1, $\{ (\tilde{x},\tilde{t})\big|\chi(\tilde{x},\tilde{t}) = -1\}$ represents fluid 2, and $\{ (\tilde{x},\tilde{t})\big|-1<\chi(\tilde{x},\tilde{t}) <1\}$  represents the immiscible two-phase flow interface.

There have been many important results about the well-posedness and dynamic behaviors of the solutions to the compressible NSAC system.
For the isentropic flow, we refer to Feireisl-Petzeltov\'{a}-Rocca-Schimperna \cite{FPRS2010} and Chen-Wen-Zhu \cite{cwz2019} for the global existence of the weak solution in three-dimensional bounded domain and \cite{chs2021-nsac,zhao2019}
for the global well-posedness and time-decay rates of the strong solutions to Cauchy problem in three dimensions. In one dimension, the existence of global classical solutions and weak solutions were investigated in \cite{DLL2013,CG2017} for the initial boundary value problems and in \cite{dlt2019,dlw2024,lydc2023} for the free boundary problems, and the stability of rarefaction waves and other wave behaviors were studied in \cite{yz2019,lyz2018,f2014}. For the non-isentropic flow, the existence and uniqueness of local strong solutions to the multi-dimensional initial boundary value problem was proved in \cite{kotschote2012} and optimal time decay rates of the global solutions to the three-dimensional Cauchy problem was shown in \cite{clt2021-non-nsac}. Concerning with the one-dimensional case, we refer to \cite{chhs2021-1d,ydl2022} for the global existence and \cite{lyz2020,luoting2022} for the stability of rarefaction wave and the composite wave. For more work on the well-posedness of compressible NSAC system, see references \cite{cz2021,f2014,kotschote2017,lc2024}, etc.

In this paper, we mainly study the large time behavior of the solutions to the initial value problem for the compressible NSAC system  with the immiscible two-phase flow initially located near the phase separation state in one dimension.
For the convenience of analysis, we introduce Lagrange coordinate transformation as follows
$$
x= \int_{(0,0)}^{(\tilde{x}, \tilde{t})}  \rho(y, s) d y-( \rho u)(y, s) d s, \quad t = \tilde{t} .
$$
Moreover, according to the far-field states \eqref{inf-e}, asymptotic state of the phase field may be $1$ or $-1$. Consequently, we redefine the new variables
\begin{equation}\label{Variable substitution 1}
  v=\rho^{-1}, \qquad \varphi=\chi^2,
\end{equation}
where $v$ is specific volume, then the Cauchy problem \eqref{nsac-e}-\eqref{init-e} can be transformed as follows in Lagrangian coordinates:
\begin{equation}\label{nsac-l}
\left\{\begin{array}{llll}
\displaystyle v_{t}-u_{x}=0,\\
\displaystyle u_{t}+p(v)_{x}=\Big(\nu(v)u_{x}\Big)_{x}-\frac{\epsilon}{8}\Big(\frac{\varphi_{x}^{2}}{\varphi v^{2}}\Big)_{x}, \\
\displaystyle \varphi_{t}=-\frac{2v}{\epsilon}(\varphi-1)\varphi+\epsilon v\Big(\frac{\varphi_{x}}{v}\Big)_{x}-\frac{\epsilon\varphi_{x}^{2}}{2\varphi},
\end{array}\right.
\end{equation}
subject to the initial data
\begin{equation}\label{init-l}
 (v,u,\varphi)(x,0)=(v_0,u_0,\varphi_0)(x)\xrightarrow{|x|\rightarrow+\infty}( \bar{v}, \bar{u}, 1),
\end{equation}
where  $v_0=\rho_0^{-1}, \varphi_0=\chi_0^2, \bar{v}=\bar{\rho}^{-1}$, and the viscosity coefficient $\nu(v)=\tilde{\nu}(v)v^{-1}>0$ is a smooth function of $v$.  The pressure $p(v)=P(\rho)$ is assumed to be a smooth function of $v$ satisfying
\begin{equation}
  p'(v)<0,\quad\text{and}\quad p''(v)>0.
\end{equation}

Our aim is to obtain a clear picture of the large-time behavior of Cauchy problem for the one-dimensional compressible NSAC system when the initial data $(v_0,u_0,\varphi_0)(x)$ is smooth and small perturbation of the constant state $(\bar{v}, \bar{u}, 1)$. 
Motivated by \cite{liu1991,zeng1994}, we introduce the following modified parabolic system corresponding to \eqref{nsac-l}-\eqref{init-l}:
\begin{equation}\label{ps}
\left\{\begin{array}{llll}
\displaystyle \tilde{v}_{t}-\tilde{u}_{x}=\frac{\nu(\bar{v})}{2}\tilde{v}_{xx},\\
\displaystyle \tilde{u}_{t}+p(\tilde{v})_{x}=\frac{\nu(\bar{v})}{2}\tilde{u}_{xx}, \\
\displaystyle \tilde{\varphi}_{t}+\epsilon (\tilde{\varphi}-1)=\epsilon\tilde{\varphi}_{xx},
\end{array}\right.
\end{equation}
with the same initial data, that is
\begin{equation}\label{init-ps}
(\tilde{v},\tilde{u},\tilde{\varphi})(x,0)=(v_0,u_0,\varphi_0)(x).
\end{equation}
It is well known that $(\tilde{v},\tilde{u})$ is time asymptotically approximated in $L^p (1 \leq p \leq+\infty)$ by diffusion waves, which were constructed in terms of the self-similar solutions to the generalized viscous Burgers equations, see \cite{cl1987,lz1997}. We also refer to \cite{k1987} for the similar results in $L^2$-norm of solutions to general hyperbolic-parabolic systems and \cite{zeng1994} for the large time behavior in $L^1$ of solutions to the compressible isentropic viscous flow. 
In this paper, we show that $(v,u)$ is time asymptotically approximated by $(\tilde{v},\tilde{u})$ at a higher rate, hence that $(v,u)$ is approximated by the same diffusion waves in \cite{zeng1994}. Moreover, we prove that $\varphi-1$ is approximated exponentially by a diffusion wave, which is the solution of heat equation with damping term (see \eqref{phi-ps}).\\

\noindent\textbf{\normalsize Notations.} 
In this paper, the operator $D^l\triangleq\frac{\partial^l}{\partial x^l}$ with an integer $l\geq 0$. For function spaces, $L^r$ denotes the space of measurable functions on $\mathds{R}$ with the norm $\|f\|_{L^r}=(\int_{\mathds{R}}|f|^rdx)^{\frac{1}{r}},\ 1\leq r<\infty$. $W^{s,r}$ denotes the Sobolev space on $\mathds{R}$ with the norm $\|f\|_{W^{s,r}}=(\sum_{l=0}^s\|{D^l f}\|_{L^r}^r)^{\frac{1}{r}}$, and $H^s \triangleq W^{s,2}$. Throughout this paper, $C$ denotes the generic positive constant and $C_i~(i=1,2,\cdots,7)$ will also denote some positive constants depending only on the given problems. We will employ the notation $A\lesssim B$ to mean that $A\leq C B$. And $A\sim B$ shows that $C^{-1}B\leq A\leq C B$.\\

To begin with, we establish the global existence and large-time behavior of the classical solutions to the system \eqref{nsac-l}-\eqref{init-l} in the following theorem:
\begin{theorem}\label{thm1}
Suppose that the initial data $(v_0,u_0,\varphi_0)$ satisfies
 $(v_0-\bar{v}, u_0-\bar{u},\varphi_0-1)\in H^s$ for an integer $s\geq 3$, and there exists a sufficiently small positive constant $\delta_0$ such that if
\begin{equation}\label{ini-1}
  \|(v_0-\bar{v}, u_0-\bar{u}, \varphi_0-1)\|_{H^s}\leq \delta_0,
\end{equation} 
then the Cauchy problem \eqref{nsac-l}-\eqref{init-l} admits a unique solution  $(v,u,\varphi)(x,t)$ satisfying
\begin{align*}
  v-\bar{v}\in C([0,+\infty),H^{s})\cap C^1([0,+\infty),H^{s-1}),\\
  u-\bar{u}\in C([0,+\infty),H^{s})\cap C^1([0,+\infty),H^{s-2}),\\
  \varphi-1\in C([0,+\infty),H^{s})\cap C^1([0,+\infty),H^{s-2}),
\end{align*}
and for any $T>0$,
\begin{equation}\label{con1}
    \sup_{0\leq t\leq T}\|(v-\bar{v}, u-\bar{u}, \varphi-1)(t)\|_{H^s}+\int_0^T(\| v_x (t)\|_{H^{s-1}}^2+ \| u_x (t)\|_{H^{s}}^2+\|(\varphi-1)(t)\|_{H^{s+1}}^2)dt \lesssim\delta_0.
\end{equation}
Assume further 
\begin{equation}\label{ini-2}
  \|(v_0-\bar{v}, u_0-\bar{u})\|_{L^1}\leq \delta_0,
\end{equation}
then for  $0\leq l\leq s$,
\begin{align}
  &  \|D^l(v-\bar{v}, u-\bar{u})(t)\|_{H^{s-l}}\lesssim\delta_0(1+t)^{-\frac12(l+\frac12)},\qquad \|(\varphi-1)(t)\|_{H^s}\lesssim\delta_0e^{-\frac{2\bar{v}t}{\epsilon}}.\label{lt1}
\end{align}
\end{theorem}

\begin{rem}\label{rem:1}
It should be noted that the time-decay rates \eqref{lt1} imply that specific volume and velocity of 1-D compressible NSAC system decays at the same time-decay rate $(1+t)^{-\frac{1}4}$ in $L^2$-norm as the compressible isentropic viscous 1-D flow in \cite{k1987} (also see \cite{zeng1994}). Furthermore, we give the time decay rates of the highest order spatial derivatives of the specific volume and velocity by utilizing the damping structure of the perturbation equation of phase field and combining the energy method with a delicate spectral analysis.
In addition, we apply the pure energy method instead of spectrum analysis to obtain the exponential decay rate of the phase field, which does not need the assumption on $L^1$-norm of the phase field. It is also shown that if the initial data is a small perturbation of the constant state, the phase field jumps and phase separation occurs for every $t>0$ as the thickness of the interface tends to zero.
\end{rem}
\begin{rem}\label{rem:2}
For the modified parabolic system \eqref{ps}-\eqref{init-ps}, we have the following results similar as that in \cite{k1987,zeng1994}. Assume that the initial data $(v_0,u_0,\varphi_0)$ satisfies \eqref{ini-1} and \eqref{ini-2}, then Cauchy problem \eqref{ps}-\eqref{init-ps} admits a unique global-in-time solution $(\tilde{v},\tilde{u},\tilde{\varphi})(x,t)$  satisfying
\begin{equation}\label{con2}
    \sup_{0\leq t\leq T}\|(\tilde{v}-\bar{v}, \tilde{u}-\bar{u}, \tilde{\varphi}-1)(t)\|_{H^s}+\int_0^T(\|  \tilde{v}_x \|_{H^{s-1}}^2+ \| \tilde{u}_x \|_{H^{s}}^2+\|\tilde{\varphi}-1\|_{H^{s+1}}^2)dt\lesssim\delta_0,
\end{equation}
and it holds that for  $0\leq l\leq s$,
\begin{align}
  &  \|D^l(\tilde{v}-\bar{v}, \tilde{u}-\bar{u})(t)\|_{H^{s-l}}\lesssim\delta_0(1+t)^{-\frac12(l+\frac12)},~ \|(\tilde{\varphi}-1)(t)\|_{H^s}\lesssim\delta_0e^{-\epsilon t}.\label{lt2}
\end{align}
\end{rem}

\

Motivated by the works in \cite{k1987,zeng1994}, we have the following large-time behavior for the difference of the solutions to the problems \eqref{nsac-l}-\eqref{init-l} and \eqref{ps}-\eqref{init-ps}.
\begin{theorem}\label{thm2}
Under the assumption of Theorem \ref{thm1}, for $\alpha>0$ arbitrarily small, $\bar{\eta}\triangleq \min\{\frac{2\bar{v}}{\epsilon}-\alpha,\epsilon\}$, and $0\leq k\leq s-1$, it holds
\begin{equation}\label{lt3}
   \|D^k(v-\tilde{v}, u-\tilde{u})(t)\|_{H^{s-1-k}}\lesssim\delta_0(1+t)^{-\frac12(k+\frac32)+\alpha},~ \|(\varphi-\tilde{\varphi})(t)\|_{H^s}\lesssim\delta_0e^{-\bar{\eta}t }.
\end{equation}
Further suppose that $(v_0-\bar{v}, u_0-\bar{u}, \varphi_0-1)\in W^{s,1} \cap H^s$ for $s\geq 3$, then we have
\begin{equation}\label{lt4}
\begin{aligned}
 &\|D^k(v-\tilde{v}, u-\tilde{u})(t)\|_{W^{s\!-\!1\!-\!k,1}}\lesssim\delta_0(1+t)^{-\frac12(k+1)+\alpha}+e^{-\frac{\bar{c}^2}{\bar{\nu}}t}\|(v_0-\bar{v}, u_0-\bar{u})\|_{W^{k,1}},\\
 & \|(\varphi-\tilde{\varphi})(t)\|_{W^{s,1}}\lesssim(\delta_0^2+\|\phi_0\|_{L^1})e^{-\bar{\eta}t}.  
\end{aligned}  
\end{equation}
\end{theorem}

\begin{rem}\label{rem-lp}
Compared with \cite{k1987,zeng1994} for the compressible viscous flow, we obtain the time decay rates of the higher order spatial derivatives $D^k(v-\tilde{v}, u-\tilde{u})$ with $0\leq k\leq s-1$.
Using the interpolation formulas, \eqref{lt3} and \eqref{lt4},
we have the following $L^p~(1\leq p\leq +\infty)$ decay rates:
\begin{equation}\label{lt5}
\begin{aligned}
 &\|D^k(v-\tilde{v}, u-\tilde{u})(t)\|_{L^p}\lesssim (1+t)^{-\frac12(k+2-\frac1p)+\alpha},\\
 & \|D^k(\varphi-\tilde{\varphi})(t)\|_{L^p}\lesssim e^{-\bar{\eta}t},  
\end{aligned}  
\end{equation}
where $0\leq k\leq s-1$. 
It is worth pointing out that $\tilde{\varphi}-1$ is a diffusion wave determined by the heat equation with damping term. \eqref{lt5}$_2$ provides that $\varphi-1$ is time asymptotically approximated in $L^p (1 \leq p \leq+\infty)$ by such diffusion wave at an exponential rate.
\end{rem}

\begin{rem}\label{rem-3}
We review that $(\tilde{v},\tilde{u})$ is time asymptotically approximated in $L^p (1 \leq p \leq+\infty)$ by diffusion waves $\Psi(x,t)=(\psi_1,\psi_2)$ in \cite{lz1997}: Assume that $(v_0,u_0)$ has the properly that
\begin{equation}\label{ini-3}
  \|(v_0-\bar{v}, u_0-\bar{u})\|_{L^1}+\|(x(v_0-\bar{v}), x(u_0-\bar{u}))\|_{L^1}+\|(v_0-\bar{v}, u_0-\bar{u})\|_{H^s}\leq \delta_0,
\end{equation} 
then it holds for $1\leq p\leq +\infty$ and $t>0$,
\begin{equation}\label{lt6}
\begin{aligned}
 &\|D^k(\tilde{v}-\psi_1, \tilde{u}-\psi_2)(t)\|_{L^p}\lesssim \delta_0t^{-\frac12(k+\frac32-\frac1p)},
\end{aligned}  
\end{equation}
where $0\leq k\leq s-1$. 
\eqref{lt5} and \eqref{lt6} immediately show that the solution $(v,u)$ of compressible NSAC system is approximated by the same diffusion waves $(\psi_1,\psi_2)$. That is to say,  under assumption \eqref{ini-3}, there are the following results:
\begin{equation}\label{lt7}
\begin{aligned}
 &\|D^k(v-\psi_1, u-\psi_2)(t)\|_{L^p}\lesssim \delta_0t^{-\frac12(k+\frac32-\frac1p)},
\end{aligned}  
\end{equation}
where $0\leq k\leq s-1$. 

\end{rem}

The rest of the paper is organized as follows.
In Section \ref{se2}, we reformulate the original problem in
terms of the perturbed variables and derive the uniform a priori estimates, which can guarantee the local classical solution to be a global classical one.
Section \ref{Sec4} is devoted to establish the decay estimates in $L^2$-norm and $L^1$-norm, respectively, and complete the proof of our main theorems.

\section{Global existence}\label{se2}
\subsection{Reformulation of the problems}
\indent\qquad
In this subsection, we reformulate the problems \eqref{nsac-l}-\eqref{init-l} and \eqref{ps}-\eqref{init-ps} in terms of the perturbated variables.  To begin with, let us define the perturbated variables as
\begin{equation}\label{pertub}
n=v-\bar{v},~ w=u-\bar{u},~ \phi=\varphi-1,
\end{equation}
\begin{equation}\label{pertub2}
\tilde{n}=\tilde{v}-\bar{v},~ \tilde{w}=\tilde{u}-\bar{u},~ \tilde{\phi}=\tilde{\varphi}-1.
\end{equation}
and introduce the constants
\begin{equation}\label{new-const}
\bar{c}=\sqrt{-p'(\bar{v})},~ \bar{\nu}=\nu(\bar{v}).
\end{equation}
The problems \eqref{nsac-l}-\eqref{init-l} and \eqref{ps}-\eqref{init-ps} are then reformulated into the following systems for the perturbations:
\begin{align}
&\left\{\begin{array}{llll}
\displaystyle n_{t}-w_{x}=0,\label{nsac-s}\\
\displaystyle w_{t}-\bar{c}^2n_{x}=\bar{\nu}w_{xx}+(f_1(n)+f_2(n,w)+f_3(n,\phi))_x, \\
\displaystyle \phi_{t}+\frac{2\bar{v}}{\epsilon}\phi=\epsilon\phi_{xx}+f_4(n,\phi)+f_5(\phi),\\
(n,w,\phi)(x,0)=(n_0,w_0,\phi_0)(x),
\end{array}\right.\\
&\left\{\begin{array}{llll}
\displaystyle \tilde{n}_{t}-\tilde{w}_{x}=\frac{\bar{\nu}}{2}\tilde{n}_{xx},\label{ps-s} \vspace{1ex}\\
\displaystyle \tilde{w}_{t}- \bar{c}^2\tilde{n}_{x}=\frac{\bar{\nu}}{2}\tilde{w}_{xx}+(\tilde{f_1}(\tilde{n}))_x,  \vspace{1ex} \\
\displaystyle \tilde{\phi}_{t}+\epsilon \tilde{\phi}=\epsilon\tilde{\phi}_{xx},\\
(\tilde{n},\tilde{w},\tilde{\phi})(x,0)=(n_0,w_0,\phi_0)(x).
\end{array}\right.
\end{align}
where the initial data $(n_0,w_0,\phi_0)(x)\triangleq(v_0(x)-\bar{v},u_0(x)-\bar{u},\varphi_0(x)-1)$ and the nonlinear term $f_{i}~(i=1,\cdots,5)$ satisfies
\begin{equation}\label{nonl}
\begin{aligned}
&f_1(n) \triangleq -p(n+\bar{v})+p(\bar{v})+p'(\bar{v})n,\\
&f_2(n,w)\triangleq (\nu(n+\bar{v})-\nu(\bar{v}))w_x,\\
&f_3(n,\phi)\triangleq -\frac{\epsilon\phi_x^2}{8(\phi+1)(n+\bar{v})^2},\\
&f_4(n,\phi)\triangleq -\frac{\epsilon n_x\phi_x}{n+\bar{v}}-\frac{2n(\phi^2+\phi)}{\epsilon},\\
&f_5(\phi)\triangleq -\frac{\epsilon\phi_x^2}{2(\phi+1)}-\frac{2\bar{v}\phi^2}{\epsilon},
\end{aligned}
\end{equation}
and $\tilde{f_1}(\tilde{n})\triangleq -p(\tilde{n}+\bar{v})+p(\bar{v})+p'(\bar{v})\tilde{n}.$

Since the local existence for the systems \eqref{nsac-s} and \eqref{ps-s} 
can be proved by the standard contraction mapping principle (see, for example, \cite{MN1980,kotschote2012}),  it suffices to obtain the following a priori estimates.
\begin{prop}\label{pro-est}
Under the assumption of Theorem \ref{thm1}, let $(n,w,\phi)$ be the classical solution of the system \eqref{nsac-s}  for $T>0$ and satisfying 
\begin{equation}\label{a-priori est}
  \sup_{0\leq t\leq T}\|(n(t),w(t),\phi(t))\|_{H^s}\leq \delta,
\end{equation}
where $s\geq 3$ is an integer and $\delta$ is a suitably small positive constant,
then it holds
  \begin{equation}\label{pro-est1}
    \|(n,w,\phi)(t)\|_{H^s}^2+\int_0^t\left(\|  n_x \|_{H^{s-1}}^2+\|  w_x \|_{H^{s}}^2+\|\phi\|_{H^{s+1}}^2\right)d\tau\lesssim \delta_0^2,
  \end{equation}  
where $\delta_0$ is a small positive constant in \eqref{ini-2}.  
\end{prop}

The proof of Proposition \ref{pro-est} follows from a series of lemmas in the next subsection. To this end, by a-priori assumption \eqref{a-priori est}, Sobolev's inequalities and the smallness of $\delta$, it is easy to get the following estimates
\begin{equation}\label{l-inf-est}
\begin{aligned}
&\frac{1}{2}\bar{v}\leq v(x,t)\leq 2\bar{v},\quad  \frac{1}{2}\leq \varphi(x,t)\leq 2,\\
&\| D^k (n,w,\phi)\|_{L^\infty}\lesssim \delta \quad \mathrm{for} ~0\leq k\leq s-1.
\end{aligned}
\end{equation}

\subsection{The a priori estimates}
\indent\qquad
In this subsection, we will establish the time-independent a priori bounds of the solutions of the reformulated system \eqref{nsac-s}. We first give the basic $L^2$ energy estimate of the solution.
\begin{lem}\label{lem-b}
Assume the conditions in Proposition \ref{pro-est} hold, then for $0<t<T$, it holds
\begin{equation}\label{basic2}
    \begin{aligned}
      &\frac{d}{dt}(\|n\|_{L^2}^2+\|w\|_{L^2}^2+\|\phi\|_{L^2}^2)
      +\bar{\nu}\|w_x\|_{L^2}^2+\frac{2\bar{v}}{\epsilon}\|\phi\|_{L^2}^2+\epsilon\| \phi_x\|_{L^2}^2\\      
      &\lesssim \delta(\|n_x\|_{L^2}^2+\|w_x\|_{L^2}^2+\|\phi\|_{L^2}^2).
    \end{aligned}
\end{equation}
\end{lem}
\begin{proof}
Multiplying $\eqref{nsac-s}_{1,2,3}$ by $\bar{c}^2n,w$ and $\phi$ respectively, integrating over $\mathds{R}$ and summing up, we obtain
\begin{equation}\label{basic1}
    \begin{aligned}
      &\frac{1}{2}\frac{d}{dt}(\bar{c}^2\|n\|_{L^2}^2+\|w\|_{L^2}^2+\|\phi\|_{L^2}^2)
      +\bar{\nu}\|w_x\|_{L^2}^2+\frac{2\bar{v}}{\epsilon}\|\phi\|_{L^2}^2+\epsilon\| \phi_x\|_{L^2}^2\\
      &=-\int (f_1+f_2+f_3) w_xdx+\int (f_4+f_5)\phi dx\\
      &\leq -\int f_1 w_xdx+C(\|n\|_{L^\infty}\|w_x\|_{L^2}^2+\|\phi_x\|_{L^\infty}\|\phi_x\|_{L^2}\|w_x\|_{L^2})\\
      &\quad +C(\|\phi\|_{L^\infty}\|n_x\|_{L^2}\|\phi_x\|_{L^2}+\|n\|_{L^\infty}\|\phi\|_{L^2}^2)\\
      &\leq -\int f_1 n_tdx+\frac{\bar{\nu}}{2}\| w_x\|_{L^2}^2+\frac{\epsilon}{2}\| \phi_x\|_{L^2}^2+\frac{\bar{v}}{\epsilon}\|\phi\|_{L^2}^2\\
      &\quad +C\delta(\|n_x\|_{L^2}^2+\|w_x\|_{L^2}^2+\|\phi\|_{L^2}^2).
    \end{aligned}
  \end{equation}
Defining
\begin{equation}\label{f1n}
  A(n)=p(\bar{v})n-\int_{\bar{v}}^{\bar{v}+n}p(s)ds,
\end{equation}
and using the Taylor expansion and \eqref{l-inf-est}, we have
\begin{equation}\label{add0}
   A(n)\sim n^2.
\end{equation}
For the first term on right side of \eqref{basic1}, it holds
\begin{equation}\label{add1}
    \begin{aligned}
  \int f_1 n_tdx&=\int(p(\bar{v})-p(\bar{v}+n)) n_tdx+\int p'(\bar{v})nn_tdx\\
  &=\frac{d}{dt}\left(\int A(n)dx-\frac{\bar{c}^2}{2}\int n^2dx\right).
  \end{aligned}
\end{equation}
Substituting \eqref{add1} into \eqref{basic1} 
yields \eqref{basic2} and completes the proof of Lemma \ref{lem-b}.
\end{proof}
\vspace{2mm}

Next, we strive to obtain the estimates for the derivatives of the solutions.
  Applying the operator $ D^l $ to the system $\eqref{nsac-s}$ with $1\leq l \leq s$, we have
\begin{equation}\label{nsac-ds}
\left\{\begin{array}{llll}
\displaystyle D^{l }n_{t}-D^{l }w_{x}=0,\\
\displaystyle D^{l }w_{t}-\bar{c}^2D^{l }n_{x}=\bar{\nu}D^{l }w_{xx}+D^{l+1}(f_1+f_2+f_3), \\
\displaystyle D^{l }\phi_{t}+\frac{2\bar{v}}{\epsilon}D^{l }\phi=\epsilon D^{l }\phi_{xx}+D^{l }f_4+D^{l }f_5.
\end{array}\right.
\end{equation}

\begin{lem}\label{lem-w}
  Assume the conditions in Proposition \ref{pro-est} hold, then for $0<t<T$ and $1\leq l \leq s$, it holds
\begin{equation}\label{lem-w1}
  \begin{aligned}
    &\frac{d}{dt}\left(\bar{c}^2\| D^l  n\|_{L^2}^2+\| D^l  w\|_{L^2}^2+\| D^l  \phi\|_{L^2}^2\right)+\bar{\nu}\| D^{l +1}  w_x \|_{L^2}^2+\frac{2\bar{v}}{\epsilon}\| D^{l } \phi\|_{L^2}^2+\epsilon\| D^{l +1} \phi\|_{L^2}^2\\    
    &\lesssim \delta\left(\| D^l  n\|_{L^2}^2+\| D^{l +1} w\|_{L^2}^2+\| D^{l +1}\phi\|_{L^2}^2+\| D^{l}\phi\|_{L^2}^2\right).  
  \end{aligned}
\end{equation}
\end{lem}
\begin{proof}
Multiplying $\eqref{nsac-ds}_{1,2}$ by $\bar{c}^2D^l  n,D^l  w$ respectively, integrating over $\mathds{R}$ and summing up, we obtain
  \begin{equation}\label{w1}
    \begin{aligned}
      &\frac{1}{2}\frac{d}{dt}(\bar{c}^2\| D^l  n\|_{L^2}^2+\| D^l  w\|_{L^2}^2)
      +\bar{\nu}\| D^{l +1} w\|_{L^2}^2\\
      &=\int D^{l+1}(f_1+f_2+f_3)  D^l  wdx\\
      &=-\int D^{l }(f_1+f_2+f_3) D^{l +1} wdx\\
      &\leq C (\|n\|_{L^\infty}\| D^{l }n\|_{L^2}+\|n\|_{L^\infty}\| D^{l }w_x\|_{L^2}+\|w_x\|_{L^\infty}\| D^{l }n\|_{L^2})\|D^{l +1} w\|_{L^2}\\
      &\quad+C(\|\phi_x\|_{L^\infty}\| D^{l }\phi_x\|_{L^2}+\|\phi_x\|_{L^\infty}\| D^{l }n\|_{L^2}+\|\phi_x\|_{L^\infty}\| D^{l }\phi\|_{L^2}) \|D^{l +1} w\|_{L^2}\\
      &\leq \frac{\bar{\nu}}{2}\| D^{l +1} w\|_{L^2}^2+C\delta(\| D^l  n\|_{L^2}^2+\| D^{l +1} w\|_{L^2}^2+\| D^{l +1}\phi\|_{L^2}^2+\| D^{l }\phi\|_{L^2}^2).
    \end{aligned}
  \end{equation}
Then we have 
\begin{equation}\label{w2}
  \begin{aligned}
    &\frac{d}{dt}(\bar{c}^2\| D^l  n\|_{L^2}^2+\| D^l  w\|_{L^2}^2)+\bar{\nu}\| D^{l +1}  w \|_{L^2}^2\\    
    &\lesssim \delta(\| D^l  n\|_{L^2}^2+\| D^{l +1} w\|_{L^2}^2+\| D^{l +1}\phi\|_{L^2}^2+\| D^{l }\phi\|_{L^2}^2).
  \end{aligned}
\end{equation}
Multiplying $\eqref{nsac-ds}_{3}$ by $D^l  \phi$, and integrating over $\mathds{R}$, we obtain
  \begin{equation}\label{ph1}
    \begin{aligned}
      &\frac{1}{2}\frac{d}{dt}(\| D^l  \phi\|_{L^2}^2)
      +\frac{2\bar{v}}{\epsilon}\| D^{l } \phi\|_{L^2}^2+\epsilon\| D^{l +1} \phi\|_{L^2}^2\\
      &=\int D^{l }f_4 D^l  \phi dx+\int D^{l }f_5 D^l  \phi dx.
    \end{aligned}
  \end{equation}
 By \eqref{a-priori est}, one have
 \begin{equation}\label{ph1-1}
    \begin{aligned}
      &\int D^{l }f_4 D^l  \phi dx
      =\int D^{l -1}(\frac{\epsilon\phi_xn_x}{n+\bar{v}}) D^{l +1}\phi dx-\int D^{l }(\frac{2n(\phi^2+\phi)}{\epsilon}) D^{l }\phi dx\\
      &\leq C (\|n_x\|_{L^\infty}\| D^{l -1}\phi_x\|_{L^2}+\|\phi_x\|_{L^\infty}\| D^{l -1}n_x\|_{L^2})\|D^{l +1} \phi\|_{L^2}\\
      &\quad+C(\|n\|_{L^\infty} \| D^{l }\phi\|_{L^2}+\|\phi\|_{L^\infty}\| D^{l }n\|_{L^2}) \|D^{l } \phi\|_{L^2}\\
      &\leq \frac{\epsilon}{4}\| D^{l +1} \phi\|_{L^2}^2+\frac{\bar{v}}{2\epsilon}\| D^{l } \phi\|_{L^2}^2+C\delta(\| D^l  n\|_{L^2}^2+\| D^{l }\phi\|_{L^2}^2).
    \end{aligned}
  \end{equation}
and
\begin{equation}\label{ph1-2}
    \begin{aligned}
      &\int D^{l }f_5 D^l  \phi dx
      =-\int D^{l }(\frac{\epsilon\phi_x^2}{2(\phi+1)} +\frac{2\bar{v}\phi^2}{\epsilon}) D^{l }\phi dx\\
      &\leq C (\|\phi_x\|_{L^\infty}\| D^{l }\phi_x\|_{L^2}+\|\phi_x\|_{L^\infty}\| D^{l }\phi\|_{L^2}+\|\phi\|_{L^\infty}\| D^{l }\phi\|_{L^2})\|D^{l } \phi\|_{L^2}\\
      &\leq \frac{\epsilon}{4}\| D^{l +1} \phi\|_{L^2}^2+\frac{\bar{v}}{2\epsilon}\| D^{l } \phi\|_{L^2}^2+C\delta\| D^{l }\phi\|_{L^2}^2.
    \end{aligned}
  \end{equation}
Substituting \eqref{ph1-1} and \eqref{ph1-2} into \eqref{ph1} yields
  \begin{equation}\label{ph2}
    \begin{aligned}
      &\frac{d}{dt}(\| D^l  \phi\|_{L^2}^2)
      +\frac{2\bar{v}}{\epsilon}\| D^{l } \phi\|_{L^2}^2+\epsilon\| D^{l +1} \phi\|_{L^2}^2
      \lesssim \delta(\| D^l  n\|_{L^2}^2+\| D^{l }\phi\|_{L^2}^2).
    \end{aligned}
  \end{equation}
Combining \eqref{w2} with \eqref{ph2} yields \eqref{lem-w} and finishes the proof of Lemma \ref{lem-w}.
\end{proof}

\begin{lem}\label{lem-n}
Assume the conditions in Proposition \ref{pro-est} hold, then there exists a positive constant $C_1>0$ so that for $0<t<T$ and $1\leq l\leq s$, it holds
\begin{equation}\label{lem-n1}
\begin{aligned}
      &\frac{\bar{c}^2}{2}\| D^l  n\|_{L^2}^2-\frac{d}{dt}\int ( D^{l -1}w\cdot D^l  n)dx\\
      &\lesssim C_1(\| D^{l +1} w\|_{L^2}^2+\| D^l  w\|_{L^2}^2)+\delta(\| D^l  n\|_{L^2}^2+\| D^{l +1}\phi\|_{L^2}^2).
    \end{aligned}
\end{equation}
\end{lem}
\begin{proof}
Applying the operator $ D^{l -1}$ to $\eqref{nsac-s}_2$, multiplying the resulting equation by $ D^{l } n$, and integrating over $\mathds{R}$, we obtain
  \begin{equation}\label{n1}
    \begin{aligned}
      &\bar{c}^2\| D^l  n\|_{L^2}^2+\frac{d}{dt}\int (- D^{l -1}w\cdot D^l  n)dx\\
      &=-\int  D^{l -1} w D^l  n_tdx+\bar{\nu}\int  D^{l +1}w  D^{l }ndx+\int  D^{l }(f_1+f_2+f_3) D^{l }ndx\\
      &\leq-\int  D^{l -1} w D^l  w_xdx+\frac{\bar{c}^2}{2}\| D^{l } n\|_{L^2}^2+C\| D^{l +1} w\|_{L^2}^2+C\delta(\| D^l  n\|_{L^2}^2+\| D^{l +1}\phi\|_{L^2}^2)\\
      &\leq \| D^l  w\|_{L^2}^2+\frac{\bar{c}^2}{2}\| D^{l } n\|_{L^2}^2+C\| D^{l +1} w\|_{L^2}^2+C\delta(\| D^l  n\|_{L^2}^2+\| D^{l +1}\phi\|_{L^2}^2).
    \end{aligned}
  \end{equation}
It holds that
\begin{equation}\label{n2}
    \begin{aligned}
      &\frac{\bar{c}^2}{2}\| D^l  n\|_{L^2}^2-\frac{d}{dt}\int ( D^{l -1}w\cdot D^l  n)dx\\
      &\lesssim C_1(\| D^{l +1} w\|_{L^2}^2+\| D^l  w\|_{L^2}^2)+\delta(\| D^l  n\|_{L^2}^2+\| D^{l +1}\phi\|_{L^2}^2).
    \end{aligned}
  \end{equation}
The proof of Lemma \ref{lem-n} is completed.
\end{proof}
\vspace{2mm}

With the above preparation, now we are in a position to state the proof of Proposition \ref{pro-est}.
\vspace{2mm}

\noindent{\it\textbf{Proof of Proposition \ref{pro-est}.}\ }
By \eqref{w2} and \eqref{ph2}, summing $l $ from $1$ to $s$, and combining with \eqref{basic2}, we obtain
\begin{equation}\label{est1}
  \begin{aligned}
    &\frac{d}{dt}\left(\|n\|_{H^{s}}^2+\|w\|_{H^{s}}^2+\|\phi\|_{H^{s}}^2\right)+\bar{\nu}\|w_x\|_{H^{s}}^2+\frac{2\bar{v}}{\epsilon}\|\phi\|_{H^{s}}^2+\epsilon\| \phi_x\|_{H^{s}}^2\\
    &\lesssim \delta (\| n_x\|_{H^{s-1}}^2+\| w_x\|_{H^{s}}^2+\|\phi\|_{H^{s+1}}^2).
  \end{aligned}
\end{equation}  
Similarly, by \eqref{n2}, summing $l $ from $1$ to $s$ gives
\begin{equation}\label{est2}
  \begin{aligned}
    &\frac{\bar{c}^2}{2}\|n_x\|_{H^{s-1}}^2-\frac{d}{dt}\left(\sum_{l =1}^s\int  D^{l -1}w\cdot D^l  ndx\right)\\
    &\lesssim C_1\| w_x\|_{H^{s}}^2+\delta(\|n_x\|_{H^{s-1}}^2+\| D^{2}\phi\|_{H^{s-1}}^2).
  \end{aligned}
\end{equation}
Choosing $\beta_1$ sufficiently small, then taking $\eqref{est1}+\beta_1\times\eqref{est2}$ yields
\begin{equation}\label{est3}
  \begin{aligned}
    &\frac{d}{dt}\left(\|n\|_{H^{s}}^2+\|w\|_{H^{s}}^2+\|\phi\|_{H^{s}}^2-\sum_{l =1}^s\int  D^{l -1}w\cdot D^l  ndx\right)\\
    &+\bar{\nu}\| w_x\|_{H^{s}}^2+\frac{2\bar{v}}{\epsilon}\|\phi\|_{H^{s}}^2+\epsilon\| \phi_x\|_{H^{s}}^2+\frac{\bar{c}^2\beta_1}{2}\| n_x\|_{H^{s-1}}^2\\
    &\lesssim \beta_1C_1\| w_x\|_{H^{s}}^2+\delta(\| n_x\|_{H^{s-1}}^2+\| w_x\|_{H^{s}}^2+\|\phi\|_{H^{s+1}}^2).
  \end{aligned}
\end{equation} 
Since $\beta_1$ and $\delta$ are sufficiently small,  one has
\begin{equation}\label{est4}
  \begin{aligned}
    &\frac{d}{dt}\left(\|n\|_{H^{s}}^2+\|w\|_{H^{s}}^2+\|\phi\|_{H^{s}}^2-\beta_1\sum_{l =1}^s\int  D^{l -1}w\cdot D^l  ndx\right)\\
    &+C_2(\|  w_x \|_{H^{s}}^2+\|\phi\|_{H^{s+1}}^2+\|  n_x \|_{H^{s-1}}^2)\leq 0,
  \end{aligned}
\end{equation} 
where $C_2$ is a positive constant. We integrate the above inequality in time
\begin{equation}\label{est5}
  \begin{aligned}
    &\left(\|n\|_{H^{s}}^2+\|w\|_{H^{s}}^2+\|\phi\|_{H^{s}}^2\right)(t)
    +C_2\int_0^t(\|  w_x \|_{H^{s}}^2+\|\phi\|_{H^{s+1}}^2+\|  n_x \|_{H^{s-1}}^2)d\tau\\
    &\lesssim \delta_0+\beta_1\sum_{l =1}^s\int | D^{l -1}w  D^l  n|dx.
  \end{aligned}
\end{equation}
Then applying the smallness of $\beta_1$ gives 
\begin{equation}\label{est6}
  \begin{aligned}
    \|(n,w,\phi)(t)\|_{H^s}^2
    +\int_0^t(\|  w_x \|_{H^{s}}^2+\|\phi\|_{H^{s+1}}^2+\|  n_x \|_{H^{s-1}}^2)d\tau\lesssim \delta_0.
  \end{aligned}
\end{equation}
Thus we finish the proof of Proposition \ref{pro-est}.
\endproof

\section{Decay estimates}\label{Sec4}
\subsection{Decay estimates for linearized problem}\label{Sec4-1}
\hspace{2em}In this subsection, we consider the linearized problem of the systems \eqref{nsac-s} and \eqref{ps-s}, which is crucial to obtain the decay estimates for the nonlinear decay estimates and the difference bewteen the solutions of the systems \eqref{nsac-s} and \eqref{ps-s}. The linearized problem are written as follows
\begin{align}
&\left\{\begin{array}{llll}
\displaystyle n_{t}-w_{x}=0,\label{nsac-sl}\\
\displaystyle w_{t}-\bar{c}^2n_{x}=\bar{\nu}w_{xx}, \\
\displaystyle \phi_{t}+\frac{2\bar{v}}{\epsilon}\phi=\epsilon\phi_{xx},\\
(n,w,\phi)(x,0)=(n_0,w_0,\phi_0)(x),
\end{array}\right.\\
&\left\{\begin{array}{llll}
\displaystyle \tilde{n}_{t}-\tilde{w}_{x}=\frac{\bar{\nu}}{2}\tilde{n}_{xx},\label{ps-sl} \vspace{1ex}\\
\displaystyle \tilde{w}_{t}- \bar{c}^2\tilde{n}_{x}=\frac{\bar{\nu}}{2}\tilde{w}_{xx},  \vspace{1ex} \\
\displaystyle \tilde{\phi}_{t}+\epsilon \tilde{\phi}=\epsilon\tilde{\phi}_{xx},\\
(\tilde{n},\tilde{w},\tilde{\phi})(x,0)=(n_0,w_0,\phi_0)(x).
\end{array}\right.
\end{align}

It is easy to verify that the equations \eqref{nsac-sl}$_3$ and \eqref{ps-sl}$_3$ have the damping structure which yield the exponential decay of $\phi$ and $\tilde{\phi}$ for the linearized systems. Hence, we first consider the decay estimates of $W=(n,w)^\top$ and $\widetilde{W}=(\tilde{n},\tilde{w})^\top$. For simplicity, we rewrite the linearized systems \eqref{nsac-sl}$_{1,2}$ and \eqref{ps-sl}$_{1,2}$  into the following systems
\begin{align}
&\partial_tW=\mathcal{L}W,\quad W|_{t=0}=W_0,\label{nsac-xi-l}\\
&\partial_t\widetilde{W}=\widetilde{\mathcal{L}}\widetilde{W},\quad \widetilde{W}|_{t=0}=W_0,\label{ps-xi-l}
\end{align}
where $W_0=(n_0,w_0)^\top$ and 
\begin{equation} 
\displaystyle \mathcal{L}:=\left(\begin{array}{ccc}
0 & \partial_{x}\\
 \bar{c}^2\partial_{x} & \bar{\nu}\partial_{xx} \notag
\end{array}\right),\quad  \widetilde{\mathcal{L}}:=\left(\begin{array}{ccc}
\frac{{\bar{\nu}}}{2} \partial_{xx} & \partial_{x} \vspace{1ex} \\
 \bar{c}^2\partial_{x} &\frac{\bar{\nu}}{2}\partial_{xx} \notag
\end{array}\right).
\end{equation}
Let $G(x,t)$ and $\widetilde{G}(x,t)$ be the Green's function of \eqref{nsac-xi-l} and \eqref{ps-xi-l}, respectively. After taking Fourier transform in $x$, we have
$$\partial_t\widehat{G}(\xi,t)=\widehat{\mathcal{L}}\widehat{G}(\xi,t),\quad \partial_t\widehat{\widetilde{G}}(\xi,t)=\widehat{\widetilde{\mathcal{L}}}\widehat{\widetilde{G}}(\xi,t),
$$
where $\widehat{\mathcal{L}}=\mathscr{F}(\mathcal{L})$ and $\widehat{\widetilde{\mathcal{L}}}=\mathscr{F}(\widetilde{\mathcal{L}})$ are given by
\begin{equation}
\widehat{\mathcal{L}}=\left(\begin{array}{ccc}
 0     & i\xi \\
 i\bar{c}^2\xi & -\bar{\nu}\xi^2 
 \end{array}\right), \quad
 \widehat{\widetilde{\mathcal{L}}}=\left(\begin{array}{ccc}
 -\frac{\bar{\nu}}{2}\xi^2    & i\xi \vspace{1ex}\\
 i\bar{c}^2\xi & -\frac{\bar{\nu}}{2}\xi^2 
 \end{array}\right).
\end{equation}
It is easy to verify that $\widehat{G}(\xi,t)$ and $\widehat{\widetilde{G}}(\xi,t)$ can be expressed as
\begin{equation}\label{g-xi}
\widehat{G}(\xi,t)=\sum_{l=1}^2e^{\lambda_l(\xi)t}P_l(\xi),\quad
\widehat{\widetilde{G}}(\xi,t)=\sum_{l=1}^2e^{\widetilde{\lambda}_l(\xi)t}\widetilde{P}_l,
\end{equation}
where 
\begin{equation}\label{lam1}
 \displaystyle \lambda_{1,2}(\xi)=-\frac{\bar{\nu}}{2}\xi^2\mp\sqrt{\frac{\bar{\nu}^2}{4}\xi^4-\bar{c}^2\xi^2}=-\frac{\bar{\nu}}{2}\xi^2\mp i\xi\sqrt{\bar{c}^2-\frac{\bar{\nu}^2}{4}\xi^2},\quad\widetilde{\lambda}_{1,2}(\xi)=-\frac{\bar{\nu}}{2}\xi^2\mp i\bar{c}\xi.
 \end{equation}
and
\begin{equation}
P_{1}(\xi)=\left(\begin{array}{ccc}
 \frac{-\lambda_{2}}{\lambda_1-\lambda_2}    & \frac{i\xi}{\lambda_1-\lambda_2}\vspace{1ex}\\
 \frac{i\bar{c}^2\xi}{\lambda_1-\lambda_2} & \frac{\lambda_1}{\lambda_1-\lambda_2}
 \end{array}\right),~~
 P_{2}(\xi)=\left(\begin{array}{ccc}
 \frac{\lambda_{1}}{\lambda_1-\lambda_2}    & \frac{-i\xi}{\lambda_1-\lambda_2} \vspace{1ex}\\
 \frac{-i\bar{c}^2\xi}{\lambda_1-\lambda_2} & \frac{-\lambda_2}{\lambda_1-\lambda_2} 
 \end{array}\right),
\end{equation}
\begin{equation}\label{pp1}
\widetilde{P}_{1}=\left(\begin{array}{ccc}
 \frac{1}{2}    & -\frac{1}{2\bar{c}}\vspace{1ex}\\
 -\frac{\bar{c}}{2} & \frac{1}{2}
 \end{array}\right),~~
\widetilde{P}_{2}=\left(\begin{array}{ccc}
  \frac{1}{2}    & \frac{1}{2\bar{c}}\vspace{1ex}\\
 \frac{\bar{c}}{2} & \frac{1}{2}
 \end{array}\right).
\end{equation}

\begin{rem}\label{rem:g-ps}
Note that $\widetilde{P}_{1}$ and $\widetilde{P}_{2}$ in \eqref{pp1} are constant matrices, the inverse Fourier transform gives us that
\begin{equation}\label{g-ps}
\begin{aligned}
\widetilde{G}(x,t)&=\mathscr{F}^{-1}\left(\widehat{\widetilde{G}}(\xi,t)\right)=\mathscr{F}^{-1}\left(e^{\tilde{\lambda}_1(\xi)t}\right)\tilde{P}_1+  \mathscr{F}^{-1}\left(e^{\tilde{\lambda}_2(\xi)t}\right)\tilde{P}_2\\
&=\frac{1}{\sqrt{\bar{\nu}t}}e^{-\frac{(x+\bar{c}t)^2}{2\bar{\nu}t}}\tilde{P}_1+\frac{1}{\sqrt{\bar{\nu}t}}e^{-\frac{(x-\bar{c}t)^2}{2\bar{\nu}t}}\tilde{P}_2.
\end{aligned}
\end{equation}
\end{rem}
We state the following lemma about the asymptotic behavior of the eigenvalues (see \cite{zeng1994,k1987}). 
\begin{lem}\label{lem-lxi}
Assume the positive constant $r_0$ is sufficiently small and $R_0$ is sufficiently large.
\begin{itemize}
  \item For low frequency part $|\xi|<r_0$, the eigenvalues $\lambda_{1,2}$ are complex conlugates satisfying
  \begin{equation}\label{lem-lxi1}
\lambda_{1,2}=\mp i\bar{c}\xi-\frac{\bar{\nu}}{2}\xi^2+\mathcal{O}(\xi^3),
\end{equation}
and
\begin{equation}
P_{1}=\left(\begin{array}{ccc}
 \frac{1}{2}+\mathcal{O}(\xi)    & -\frac{1}{2\bar{c}}+\mathcal{O}(\xi^2) \vspace{1ex}\\
 -\frac{\bar{c}}{2}+\mathcal{O}(\xi^2) & \frac{1}{2}+\mathcal{O}(\xi)  
 \end{array}\right),~~
 P_{2}=\left(\begin{array}{ccc}
 \frac{1}{2}+\mathcal{O}(\xi)    & \frac{1}{2\bar{c}}+\mathcal{O}(\xi^2)  \vspace{1ex}\\
 \frac{\bar{c}}{2}+\mathcal{O}(\xi^2) & \frac{1}{2}+\mathcal{O}(\xi)  
 \end{array}\right).
\end{equation}
In addition, there exists a positive constant $\hat{\nu}$ such that
\begin{equation}\label{lem-lxi2}
Re\lambda_l\leq -\hat{\nu}|\xi|^2, \quad l=1,2.
\end{equation}
  \item For high frequency part $|\xi|>R_0$, the eigenvalues are real and satisfy
\begin{equation}\label{lem-hxi1}
\lambda_1=-\bar{\nu}\xi^{2}+\frac{\bar{c}^2}{\bar{\nu}}+\mathcal{O}(\xi^{-2}),~~\lambda_2=-\frac{\bar{c}^2}{\bar{\nu}}+\mathcal{O}(\xi^{-2}),
\end{equation}
and
\begin{equation}
P_{1}=\left(\begin{array}{ccc}
 \mathcal{O}(\xi^{-2})    &\!\! -\frac{i}{\bar{\nu}\xi}\!+\!\mathcal{O}(\xi^{-3})  \\
 -\frac{i\bar{c}^2}{\bar{\nu}\xi}\!+\!\mathcal{O}(\xi^{-2}) &\!\! 1\!+\!\mathcal{O}(\xi^{-2}) 
 \end{array}\right),~~
 P_{2}=\left(\begin{array}{ccc}
  1\!+\!\mathcal{O}(\xi^{-2})    &\!\! \frac{i}{\bar{\nu}\xi}\!+\!\mathcal{O}(\xi^{-3})  \\
 \frac{i\bar{c}^2}{\bar{\nu}\xi}\!+\!\mathcal{O}(\xi^{-3}) &\!\! \mathcal{O}(\xi^{-2})  
 \end{array}\right).
\end{equation}
In addition, there exists a positive constant $R_1=\min\{\bar{\nu}R_0^2,\frac{\bar{c}^2}{\bar{\nu}}\}$ such that
\begin{equation}\label{lem-hxi2}
Re\lambda_l\leq -R_1, \quad l=1,2.
\end{equation}
\item For medium frequency part $r_0\leq|\xi|\leq R_0$, there exists a positive constant $R_2=\min\{\bar{\nu}r_0^2,\frac{\bar{c}^2}{\bar{\nu}}\}$ such that the eigenvalues satisfy
\begin{equation}\label{lem-mxi1}
Re\lambda_l\leq -R_2,\quad l=1,2.
\end{equation}
\end{itemize}
\end{lem}

\begin{rem}\label{rem-gg}
{\rm (1)} By \eqref{g-xi}-\eqref{pp1} and the first part of Lemma \ref{lem-lxi}, we obtain that for the low frequency part $|\xi|<r_0$ and $t\geq 0$,
\begin{equation}\label{gg1}
\big|\widehat{G}(\xi,t)-\widehat{\widetilde{G}}(\xi,t)\big|\lesssim |\xi|e^{-\hat{\nu}\xi^2t}.
\end{equation}
It implies that $\widehat{G}(\xi,t)$ is well approximated by $\widehat{\widetilde{G}}(\xi,t)$ for $\xi\rightarrow 0$.

{\rm (2)} From \eqref{g-xi} and the second part of Lemma \ref{lem-lxi}, we see that as $\xi\rightarrow \infty$ along the upper edge of the real axis,
\begin{equation}\label{gg2}
  e^{\lambda_1(\xi)t}P_1(\xi)\rightarrow 0,\quad e^{\lambda_2(\xi)t}P_2(\xi)\rightarrow e^{-\frac{\bar{c}^2}{\bar{\nu}}t}A,
\end{equation}
where $A=\left(\begin{array}{ccc}
  1  & 0  \\
 0 & 0  \end{array}\right).$
The difference between $G(x,t)$ and $\widetilde{G}(x,t)$ satisfies the following
estimates (see \cite[Lemma 2.8]{zeng1994}):
\begin{equation}\label{gg3}
  |{G}(x,t)-\widetilde{G}(x,t)-e^{-\frac{\bar{c}^2}{\bar{\nu}}t}\delta(x)A|\lesssim (1+t)^{-\frac12}t^{-\frac12}\left(e^{-\frac{(x+\bar{c}t)^2}{Ct}}+e^{-\frac{(x-\bar{c}t)^2}{Ct}}\right),
\end{equation}
where $\delta(x)$ is the Dirac delta function.
\end{rem}

Now we are in the position to give the $L^2$-decay rate of the solutions to the linearized systems \eqref{nsac-sl} and \eqref{ps-sl}.
\begin{lem}\label{pro-linear}
Assume the initial data $W_0=(n_0,w_0)(x)\in H^s\cap L^1, \phi_0\in H^s$ with $s\geq 3$, the global solutions $(W,\phi)=(n,w,\phi)$ and $(\widetilde{W},\tilde{\phi})=(\tilde{n},\tilde{w},\tilde{\phi})$ of the linearized  problems \eqref{nsac-sl} and \eqref{ps-sl} satisfy that for $0\leq l\leq s$, 
  \begin{align}
  &  \| (D^lW,D^l\widetilde{W})(t)\|_{L^2}\lesssim \delta_0(1+t)^{-\frac{1}{4}-\frac{l}{2}},
  \label{pro-linear1}  \\
  &  \| D^l\phi(t)\|_{L^2}\lesssim \delta_0e^{-\frac{2\bar{v}}{\epsilon}t},\quad \| D^l\tilde{\phi}(t)\|_{L^2}\lesssim \delta_0e^{-\epsilon t}.\label{pro-linear2}
  \end{align}
Moreover, it holds
\begin{equation}\label{pro-linear3}
  \| D^l(W-\widetilde{W})(t)\|_{L^2}\lesssim \delta_0(1+t)^{-\frac{3}{4}-\frac{l}{2}},\quad \| D^l(\phi-\tilde{\phi})(t)\|_{L^2}\lesssim \delta_0e^{-\eta t},
\end{equation}
where $\eta \triangleq\min\{\frac{2\bar{v}}{\epsilon},\epsilon\}.$
  \end{lem}
\begin{proof}
Due to Lemma \ref{lem-lxi} and Parseval's equality, for given function $f(x,t)$, we have the following decay estimate
  \begin{equation}\label{pro-g0}
    \| D^l G*f\|_{L^2}\lesssim (1+t)^{-\frac{1}{4}-\frac{l}{2}}\|f\|_{L^1}+e^{-Rt}\| D^l f\|_{L^2}.
  \end{equation}
By the assumption of the initial data of the system \eqref{nsac-xi-l}, it holds
\begin{equation*}
\begin{aligned}
 \| D^l W(t)\|_{L^2}= \| D^lG*W_0\|_{L^2}&\lesssim (1+t)^{-\frac{1}{4}-\frac{l}{2}}\|W_0\|_{L^1}+e^{-Rt}\| D^lW_0\|_{L^2}\\
 &\lesssim \delta_0(1+t)^{-\frac{1}{4}-\frac{l}{2}}. 
\end{aligned}
    \end{equation*}
Similarly, by \eqref{g-ps}, \eqref{gg1} and \eqref{gg3}, one has
\begin{equation*}
\begin{aligned}
 &\| D^l \widetilde{W}(t)\|_{L^2}= \| D^l\widetilde{G}*W_0\|_{L^2}\lesssim \delta_0(1+t)^{-\frac{1}{4}-\frac{l}{2}},\\
& \| D^l (W-\widetilde{W})(t)\|_{L^2}= \| D^l(G-\widetilde{G})*W_0\|_{L^2}\lesssim \delta_0(1+t)^{-\frac{3}{4}-\frac{l}{2}}.
 \end{aligned}
\end{equation*}

It remains to establish the large time behavior of $\phi$ and $\tilde{\phi}$.
Applying the operator $ D^l$ to $\eqref{nsac-sl}_3$, multiplying the resulting equation by $ D^l\phi$, and integrating over $\mathds{R}$, one derives
  \begin{equation*}
    \frac{1}{2}\frac{d}{dt}\| D^l \phi\|_{L^2}^2+\frac{2\bar{v}}{\epsilon}\| D^l\phi\|_{L^2}^2
    +\epsilon\| D^{l+1}\phi\|_{L^2}^2=0.
  \end{equation*}
By Gr$\ddot{\text{o}}$nwall's inequality, one gets
\begin{equation}\label{phi-decay1}
  \| D^l\phi\|_{L^2}\leq e^{-\frac{2\bar{v}}{\epsilon}t}\| D^l\phi_0\|_{L^2}
  \leq e^{-\frac{2\bar{v}}{\epsilon}t}\delta_0.
\end{equation}
Similarly, from $\eqref{ps-sl}_3$, one has
\begin{equation*}
    \frac{1}{2}\frac{d}{dt}\| D^l \tilde{\phi}\|_{L^2}^2+\epsilon \| D^l\tilde{\phi}\|_{L^2}^2
    +\epsilon\| D^{l+1}\tilde{\phi}\|_{L^2}^2=0,
  \end{equation*}
which implies 
\begin{equation}\label{phi-decay2}
  \| D^l\tilde{\phi}\|_{L^2}\leq  e^{-\epsilon t}\delta_0.
\end{equation}
Finally, subtracting $\eqref{nsac-sl}_3$ from $\eqref{ps-sl}_3$ yields
\begin{equation}\label{phiphi1}
  (\phi-\tilde{\phi})_{t}+\frac{2\bar{v}}{\epsilon}\phi-\epsilon\tilde{\phi}=\epsilon(\phi-\tilde{\phi})_{xx}.
\end{equation}
Let's first consider the case $\frac{2\bar{v}}{\epsilon}\geq\epsilon$.
Applying the operator $ D^l$ to \eqref{phiphi1}, multiplying the resulting equation by $ D^l(\phi-\tilde{\phi})$, and integrating over $\mathds{R}$, one obtains
\begin{equation*}
\begin{aligned}
 &\frac{1}{2}\frac{d}{dt}\| D^l (\phi-\tilde{\phi})\|_{L^2}^2+\epsilon \| D^l(\phi-\tilde{\phi})\|_{L^2}^2
    +\epsilon\| D^{l+1}(\phi-\tilde{\phi})\|_{L^2}^2\\
 &= -(\frac{2\bar{v}}{\epsilon}-\epsilon)\int D^l\phi(D^l\phi-D^l\tilde{\phi})dx\\
 &\leq-(\frac{2\bar{v}}{\epsilon}-\epsilon)\int D^l\phi D^l\tilde{\phi}dx\\
 &\lesssim \| D^l\phi\|_{L^2}\|D^l \tilde{\phi}\|_{L^2},    
\end{aligned}
   \end{equation*}
which together with \eqref{phi-decay1} and \eqref{phi-decay2} implies that
\begin{equation*}
  \frac{d}{dt}\| D^l (\phi-\tilde{\phi})\|_{L^2}^2+2\epsilon \| D^l(\phi-\tilde{\phi})\|_{L^2}^2
    +2\epsilon\| D^{l+1}(\phi-\tilde{\phi})\|_{L^2}^2\lesssim e^{-(\frac{2\bar{v}}{\epsilon}+\epsilon)t}\delta_0^2.
\end{equation*}
Note that $\phi(x,0)=\tilde{\phi}(x,0)=\phi_0(x)$, applying Gr$\ddot{\text{o}}$nwall's inequality to the above inequality, for the case $\frac{2\bar{v}}{\epsilon}\geq\epsilon$,  one has 
\begin{equation*}
  \| D^l (\phi-\tilde{\phi})\|_{L^2}^2\lesssim \int_0^te^{-2\epsilon(t-\tau)}e^{-(\frac{2\bar{v}}{\epsilon}+\epsilon )\tau}\delta_0^2d\tau\lesssim \delta_0^2e^{-2\epsilon t},
\end{equation*}
which implies that
  $\| D^l (\phi-\tilde{\phi})\|_{L^2}\lesssim \delta_0e^{-\epsilon t}.$
Similarly, for the other case $\frac{2\bar{v}}{\epsilon}<\epsilon$, we have
  $\| D^l (\phi-\tilde{\phi})\|_{L^2}\lesssim \delta_0e^{-\frac{2\bar{v}}{\epsilon} t}.$
Then for $\eta \triangleq \min\{\frac{2\bar{v}}{\epsilon},\epsilon\}$, it holds
\begin{equation}\label{phi-decay3}
  \| D^l (\phi-\tilde{\phi})\|_{L^2}\lesssim \delta_0e^{-\eta t}.
\end{equation}
This completes the proof.
\end{proof}
\vspace{2mm}

Next, we are in the position to give the $L^1$-decay rate of the difference between the solution $(n,w,\phi)$ of the linearized system \eqref{nsac-sl} and the solution $(\tilde{n},\tilde{w},\tilde{\phi})$ of the linearized system \eqref{ps-sl}.
\begin{lem}\label{pro-l1}
Assume the initial data $(W_0,\phi_0)\in H^s\cap W^{s,1}$ with $s\geq 3$, the global solutions $(W,\phi)$ and $(\widetilde{W},\tilde{\phi})$ of the linearized  problems \eqref{nsac-sl} and \eqref{ps-sl} satisfy that for $0\leq l\leq s$, 
  \begin{align}
  &  \| D^l(W-\widetilde{W})(t)\|_{L^1}\lesssim \delta_0(1+t)^{-\frac{1}{2}-\frac{l}{2}}+e^{-\frac{\bar{c}^2}{\bar{\nu}}t}\|W_0\|_{W^{l,1}}, \label{pro-linear4}  \\
  &\| D^l(\phi-\tilde{\phi})(t)\|_{L^1}\lesssim e^{-\eta t}\|\phi_0\|_{L^1}, ~
  \| D^l\phi(t)\|_{L^1}\lesssim e^{-\frac{2\bar{v}}{\epsilon}t}\|\phi_0\|_{L^1}, ~ \| D^l\tilde{\phi}(t)\|_{L^1}\lesssim e^{-\epsilon t}\|\phi_0\|_{L^1}.\label{pro-linear5} 
  \end{align}
where $\eta \triangleq\min\{\frac{2\bar{v}}{\epsilon},\epsilon\}.$
Moreover, it holds
 \begin{align}
 &  \| (D^lW,D^l\widetilde{W})(t)\|_{L^1}\lesssim \delta_0(1+t)^{-\frac{l}{2}}+e^{-\frac{\bar{c}^2}{\bar{\nu}}t}\|W_0\|_{W^{l,1}}.\label{pro-linear6} 
\end{align}
  \end{lem}
\begin{proof}
First, we prove the estimates \eqref{pro-linear4} for $l=0$. By \eqref{gg3}, one has
\begin{equation*}
\begin{aligned}
 &\|(W-\widetilde{W})(t)\|_{L^1}= \|(G-\widetilde{G})*W_0\|_{L^1}\\
 &\lesssim \|(G-\widetilde{G}-e^{-\frac{\bar{c}^2}{\bar{\nu}}t}\delta(x)A)*W_0\|_{L^1}+e^{-\frac{\bar{c}^2}{\bar{\nu}}t}\|W_0\|_{L^1}\\
 &\lesssim \left(\|(G-\widetilde{G}-e^{-\frac{\bar{c}^2}{\bar{\nu}}t}\delta(x)A)\|_{L^1}+e^{-\frac{\bar{c}^2}{\bar{\nu}}t}\right)\|W_0\|_{L^1}\\
 &\lesssim \delta_0(1+t)^{-\frac{1}{2}}+e^{-\frac{\bar{c}^2}{\bar{\nu}}t}\|W_0\|_{L^1}.
\end{aligned}
    \end{equation*}
The difference between $G(x,t)$ and $\widetilde{G}(x,t)$ satisfies the following
estimates (see \cite[Theorem 5.8]{lz1997}): for any integer $l\geq 0$,
\begin{equation}\label{gg4}
  |D^l{G}(x,t)-D^l\widetilde{G}(x,t)-e^{-\frac{\bar{c}^2}{\bar{\nu}}t}\sum_{i=0}^l\delta^{i}(x)A_i(t)|\lesssim (1+t)^{-\frac12}t^{-\frac{1+l}{2}}\left(e^{-\frac{(x+\bar{c}t)^2}{Ct}}+e^{-\frac{(x-\bar{c}t)^2}{Ct}}\right),
\end{equation}
where $\delta^{i}(x)$ is the $i$-th derivative of the Dirac delta function and $A_i = A_i(t)$ is a $2\times2$ polynomial matrix. Additionally, we have
$A_0=\left(\begin{array}{ccc}
  1  & 0  \\
 0 & 0  \end{array}\right), A_1=\left(\begin{array}{ccc}
  0  & -\frac{1}{\bar{v}}  \\
 -\frac{\bar{c}^2}{\bar{v}} & 0  \end{array}\right).$
Hence we have
\begin{equation}\label{ww1}
\begin{aligned}
 &\|D^l(W-\widetilde{W})(t)\|_{L^1}= \|D^l(G-\widetilde{G})*W_0\|_{L^1}\\
 &\lesssim\|(D^l{G}(x,t)-D^l\widetilde{G}(x,t)-e^{-\frac{\bar{c}^2}{\bar{\nu}}t}\sum_{i=0}^l\delta^{i}(x)A_i(t))*W_0\|_{L^1}+e^{-\frac{\bar{c}^2}{\bar{\nu}}t}\sum_{i=0}^l\|A_iW_0^{(l-i)}\|_{L^1}\\
 &\lesssim \left((1+t)^{-\frac{1+l}{2}}\|W_0\|_{L^1}+e^{-\frac{\bar{c}^2}{\bar{\nu}}t}\|W_0\|_{W^{l,1}}\right)\\
 &\lesssim \delta_0(1+t)^{-\frac{1+l}{2}}+e^{-\frac{\bar{c}^2}{\bar{\nu}}t}\|W_0\|_{W^{l,1}}.
\end{aligned}
    \end{equation}
Next, by \eqref{g-ps}, a directly calculation yields
\begin{equation}\label{ww2}
\begin{aligned}
\| D^l \widetilde{W}(t)\|_{L^1}= \| D^l\widetilde{G}*W_0\|_{L^1}&\lesssim (1+t)^{-\frac{l}{2}}\|W_0\|_{L^1}++e^{-\frac{\bar{c}^2}{\bar{\nu}}t}\|W_0\|_{W^{l,1}}\\
&\lesssim \delta_0(1+t)^{-\frac{l}{2}}+e^{-\frac{\bar{c}^2}{\bar{\nu}}t}\|W_0\|_{W^{l,1}}.
 \end{aligned}
\end{equation}
The $L^1$-estimates for $D^l W$ is a consequence of \eqref{ww1} and \eqref{ww2}.

It remains to establish $L^1$-estimates for $D^l\phi$ and $D^l\tilde{\phi}$. After taking Fourier transform with respect to $x$ to the equations $\eqref{nsac-sl}_3$ and $\eqref{ps-sl}_3$, we have
\begin{align*}
\partial_t\widehat{\phi}(\xi,t)+\frac{2\bar{v}}{\epsilon}\widehat{\phi}(\xi,t)=-\epsilon\xi^2\widehat{\phi}(\xi,t),\\
\partial_t\widehat{\tilde{\phi}}(\xi,t)+\epsilon \widehat{\tilde{\phi}}(\xi,t)=-\epsilon\xi^2\widehat{\tilde{\phi}}(\xi,t).
 \end{align*}
The inverse Fourier transform gives us
\begin{align}\label{phi-ps}
\phi(x,t)=e^{-\frac{2\bar{v}}{\epsilon}t}H(x,t)*\phi_0(x),\quad
\tilde{\phi}(x,t)=e^{-\epsilon t}H(x,t)*\phi_0(x),
 \end{align}
where the heat kernel $H(x,t)=\frac{1}{\sqrt{4\pi\epsilon t}}e^{-\frac{x^2}{4\epsilon t}}$. A directly calculation yields
 \begin{align*}
 \| D^l\phi(t)\|_{L^1}\lesssim \|\phi_0\|_{L^1}e^{-\frac{2\bar{v}}{\epsilon}t}, \quad \| D^l\tilde{\phi}(t)\|_{L^1}\lesssim \|\phi_0\|_{L^1} e^{-\epsilon t}.
\end{align*}
Set $\eta \triangleq\min\{\frac{2\bar{v}}{\epsilon},\epsilon\}$, \eqref{phi-ps} implies that
\begin{align*}
  |\phi-\tilde{\phi}|=|(e^{-\frac{2\bar{v}}{\epsilon}t}-e^{-\epsilon t})H(x,t)*\phi_0(x)|
  \lesssim e^{-\eta t}|H(x,t)*\phi_0(x)|,
\end{align*}
which implies that $\| D^l(\phi-\tilde{\phi})(t)\|_{L^1}\lesssim \|\phi_0\|_{L^1}e^{-\eta t}$.
This completes the proof.
\end{proof}

\subsection{Time decay rates for nonlinear system}
\hspace{2em} In this subsection,  we are devoted to the study of the decay rates to nonlinear systems \eqref{nsac-s} and \eqref{ps-s}. 
Due to the damping structure of $\phi$, it inspires us to establish the decay of $L^2$-estimates for $\phi$ and $\phi-\tilde{\phi}$ firstly.
\begin{lem}\label{lem-pd}
Under the assumptions of Theorem \ref{thm1}, $(n,w,\phi)$ and $(\tilde{n},\tilde{\phi},\tilde{w})$ are the classical solution of the problem \eqref{nsac-s} and \eqref{ps-s}, then we have for $\alpha>0$ arbitrarily small,
  \begin{align}
   & \|\phi(t)\|_{H^s}\lesssim e^{-(\frac{2\bar{v}}{\epsilon}-\alpha)t}\delta_0,
  \quad \|(\phi-\tilde{\phi})(t)\|_{H^s}\lesssim e^{-\bar{\eta} t}\delta_0,
  \end{align}
  where $\bar{\eta}\triangleq \min\{\frac{2\bar{v}}{\epsilon}-\alpha,\epsilon\}.$
\end{lem}
\begin{proof}
To begin with, we establish the decay estimates of $\phi$. Applying the operator $D^{l }$ to $\eqref{nsac-s}_3$  and taking $L^2$ inner product with $D^l \phi$, one has 
  \begin{equation}\label{phi-1}
    \begin{aligned}
    &\frac{1}{2}\frac{d}{dt}\|D^l  \phi\|_{L^2}^2+\frac{2\bar{v}}{\epsilon}\|D^l \phi \|_{L^2}^2
    +\epsilon\|D^{l +1}\phi \|_{L^2}^2\\
    &=\int D^{l }f_4 D^l  \phi dx+\int D^{l }f_5 D^l  \phi dx.
    \end{aligned}
  \end{equation}
Via Proposition \ref{pro-est} and \eqref{l-inf-est}, the first term on the right-hand side of \eqref{phi-1} is estimated as follows
  \begin{equation*}
    \begin{aligned}
      &\int D^{l }f_4 D^l  \phi dx
      =\int D^{l -1}(\frac{\epsilon\phi_xn_x}{n+\bar{v}}) D^{l +1}\phi dx-\int D^{l }(\frac{2n(\phi^2+\phi)}{\epsilon}) D^{l }\phi dx\\
      &\lesssim (\|n_x\|_{L^\infty}\| D^{l -1}\phi_x\|_{L^2}+\|\phi_x\|_{L^\infty}\| D^{l -1}n_x\|_{L^2})\|D^{l +1} \phi\|_{L^2}\\
      &\quad+C(\|n\|_{L^\infty} \| D^{l }\phi\|_{L^2}+\|\phi\|_{L^\infty}\| D^{l }n\|_{L^2}) \|D^{l } \phi\|_{L^2}\\
      &\lesssim  \delta_0(\| D^{l +1} \phi\|_{L^2}^2+\| D^{l } \phi\|_{L^2}^2+\|\phi_x\|_{H^1}^2).
    \end{aligned}
  \end{equation*}
where we have used the fact that $\| D^l  n\|_{L^2}\lesssim \delta_0$ for $l \leq s$. Similarly, the remaining term on the right-hand side of \eqref{phi-1} are calculated as
\begin{equation*} 
\begin{aligned}
    &\int D^{l }f_5 D^l  \phi dx
      =-\int D^{l }(\frac{\epsilon\phi_x^2}{2(\phi+1)} +\frac{2\bar{v}\phi^2}{\epsilon}) D^{l }\phi dx\\
    &\lesssim (\|\phi_x\|_{L^\infty}\| D^{l }\phi_x\|_{L^2}+\|\phi_x\|_{L^\infty}\| D^{l }\phi\|_{L^2}+\|\phi\|_{L^\infty}\| D^{l }\phi\|_{L^2})\|D^{l } \phi\|_{L^2}\\
    &\lesssim  \delta_0(\| D^{l +1} \phi\|_{L^2}^2+\| D^{l } \phi\|_{L^2}^2).
    \end{aligned}
  \end{equation*}
As a consequence of the above estimates, it is easy to prove  
  \begin{equation}
    \begin{aligned}
    &\frac{1}{2}\frac{d}{dt}\|D^l  \phi\|_{L^2}^2+\frac{2\bar{v}}{\epsilon}\|D^l \phi \|_{L^2}^2+\epsilon\|D^{l +1}\phi \|_{L^2}^2\\
    &\lesssim  \delta_0(\| D^{l +1} \phi\|_{L^2}^2+\| D^{l } \phi\|_{L^2}^2+\|\phi_x\|_{H^1}^2).
     \end{aligned}  
     \end{equation}
Summing $l $ from $0$ to $s$, and noting the smallness of $\delta_0$, one has, for $\alpha>0$ arbitrarily small,
   \begin{equation}
    \frac{1}{2}\frac{d}{dt}\|\phi\|_{H^s}^2+(\frac{2\bar{v}}{\epsilon}-\alpha)\|\phi\|_{H^s}^2+\frac{\epsilon}{2}\|\phi_x\|_{H^s}^2\leq 0.  
\end{equation}
By Gr$\ddot{\text{o}}$nwall's inequality, it is easy to verify 
\begin{equation}\label{phi-decay4}
  \|\phi\|_{H^s}\lesssim e^{-(\frac{2\bar{v}}{\epsilon}-\alpha)t}\|\phi_0\|_{H^s}\lesssim e^{-(\frac{2\bar{v}}{\epsilon}-\alpha)t}\delta_0.
\end{equation}

It remains to establish the decay estimates of $\phi-\tilde{\phi}$.
Subtracting $\eqref{nsac-s}_3$ from $\eqref{ps-s}_3$ yields
\begin{equation}\label{phiphi2}
  (\phi-\tilde{\phi})_{t}+\frac{2\bar{v}}{\epsilon}\phi-\epsilon\tilde{\phi}=\epsilon(\phi-\tilde{\phi})_{xx}+f_4+f_5.
\end{equation}
Let's first consider the case $\frac{2\bar{v}}{\epsilon}-\alpha>\epsilon$.
Applying the operator $ D^l$ to \eqref{phiphi1}, multiplying the resulting equation by $ D^l(\phi-\tilde{\phi})$, and integrating over $\mathds{R}$, one obtains
\begin{equation}\label{pp-1}
\begin{aligned}
 &\frac{1}{2}\frac{d}{dt}\| D^l (\phi-\tilde{\phi})\|_{L^2}^2+\epsilon \| D^l(\phi-\tilde{\phi})\|_{L^2}^2
    +\epsilon\| D^{l+1}(\phi-\tilde{\phi})\|_{L^2}^2\\
 &= -(\frac{2\bar{v}}{\epsilon}-\epsilon)\int D^l\phi (D^l\phi-D^l\tilde{\phi})dx+\int D^l(f_4+f_5) D^l(\phi-\tilde{\phi})dx\\
 &\leq-(\frac{2\bar{v}}{\epsilon}-\epsilon)\int D^l\phi D^l\tilde{\phi}dx+\int D^lf_4 D^l(\phi-\tilde{\phi})dx+\int D^lf_5 D^l(\phi-\tilde{\phi})dx.   
\end{aligned}
   \end{equation}
Via Proposition \ref{pro-est} and \eqref{l-inf-est}, the second term on the right-hand side of \eqref{pp-1} is estimated as follows
  \begin{equation*}
    \begin{aligned}
      &\int D^{l }f_4 D^l (\phi-\tilde{\phi})dx
      =\int D^{l -1}(\frac{\epsilon\phi_xn_x}{n+\bar{v}}) D^{l +1}(\phi-\tilde{\phi}) dx-\int D^{l }(\frac{2n(\phi^2+\phi)}{\epsilon}) D^{l }(\phi-\tilde{\phi}) dx\\
      &\lesssim \delta_0(\| D^{l -1}\phi_x\|_{L^2}+\|\phi_x\|_{L^\infty})\|D^{l +1}(\phi-\tilde{\phi})\|_{L^2}+\delta_0(\| D^{l }\phi\|_{L^2}+\|\phi\|_{L^\infty})(\| D^{l }\phi\|_{L^2}+\|D^{l } \tilde{\phi}\|_{L^2})\\
      &\lesssim  \delta_0(\| D^{l +1} (\phi-\tilde{\phi})\|_{L^2}^2+\| D^{l } \phi\|_{L^2}^2+\|\phi_x\|_{H^1}^2+(\| D^{l }\phi\|_{L^2}+\|\phi\|_{H^1})\|D^{l } \tilde{\phi}\|_{L^2}).
    \end{aligned}
  \end{equation*}
Similarly, the last term on the right-hand side of \eqref{phi-1} are calculated as
\begin{equation*} 
\begin{aligned}
    &\int D^{l }f_5 D^l (\phi-\tilde{\phi})dx
      =\int D^{l-1}(\frac{\epsilon\phi_x^2}{2(\phi+1)})  D^{l+1}(\phi-\tilde{\phi}) dx-\int D^{l }(\frac{2\bar{v}\phi^2}{\epsilon}) D^{l }(\phi-\tilde{\phi}) dx\\
    &\lesssim \delta_0(\| D^{l }\phi_x\|_{L^2}\| D^{l+1}(\phi-\tilde{\phi})\|_{L^2}+\| D^{l }\phi\|_{L^2}(\| D^{l }\phi\|_{L^2}+\|D^{l } \tilde{\phi}\|_{L^2}))\\
    &\lesssim  \delta_0(\| D^{l +1} (\phi-\tilde{\phi})\|_{L^2}^2+\| D^{l } \phi\|_{L^2}^2+\| D^{l }\phi\|_{L^2}\|D^{l } \tilde{\phi}\|_{L^2}).
    \end{aligned}
  \end{equation*}
As a consequence of the above estimates, and noting the smallness of $\delta_0$, it holds
\begin{equation}
    \begin{aligned}
    &\frac{1}{2}\frac{d}{dt}\|D^l (\phi-\tilde{\phi})\|_{L^2}^2+\epsilon\|D^l (\phi-\tilde{\phi})\|_{L^2}^2+\frac{\epsilon}{2}\|D^{l +1}(\phi-\tilde{\phi}) \|_{L^2}^2\\
    &\lesssim  \delta_0(\| D^{l } \phi\|_{L^2}^2+\|\phi_x\|_{H^1}^2+(\| D^{l }\phi\|_{L^2}+\|\phi\|_{H^1})\|D^{l } \tilde{\phi}\|_{L^2})),
     \end{aligned}  
     \end{equation}
which together with \eqref{phi-decay4} and \eqref{phi-decay2} implies that
\begin{equation*}
  \frac{d}{dt}\| D^l (\phi-\tilde{\phi})\|_{L^2}^2+2\epsilon \| D^l(\phi-\tilde{\phi})\|_{L^2}^2
    +\epsilon\| D^{l+1}(\phi-\tilde{\phi})\|_{L^2}^2\lesssim (e^{-(\frac{4\bar{v}}{\epsilon}-2\alpha)t}+e^{-(\frac{2\bar{v}}{\epsilon}-\alpha+\epsilon)t})\delta_0^2.
\end{equation*}
Applying Gr$\ddot{\text{o}}$nwall's inequality to the above inequality, for the case $\frac{2\bar{v}}{\epsilon}-\alpha>\epsilon$, one has
\begin{equation*}
  \| D^l (\phi-\tilde{\phi})\|_{L^2}^2\lesssim \int_0^te^{-2\epsilon(t-\tau)}(e^{-(\frac{4\bar{v}}{\epsilon}-2\alpha)\tau}+e^{-(\frac{2\bar{v}}{\epsilon}-\alpha+\epsilon)\tau})\delta_0^2d\tau\lesssim \delta_0^2e^{-2\epsilon t},
\end{equation*}
which implies that  $\| \phi-\tilde{\phi}\|_{H^s}\lesssim \delta_0e^{-\epsilon t}.$

Next, for the other case $\frac{2\bar{v}}{\epsilon}-\alpha<\epsilon$, one similarly obtains
\begin{equation}\label{pp-2}
\begin{aligned}
 &\frac{1}{2}\frac{d}{dt}\| D^l (\phi-\tilde{\phi})\|_{L^2}^2+\frac{2\bar{v}}{\epsilon} \| D^l(\phi-\tilde{\phi})\|_{L^2}^2
    +\epsilon\| D^{l+1}(\phi-\tilde{\phi})\|_{L^2}^2\\
 &= -(\frac{2\bar{v}}{\epsilon}-\epsilon)\int D^l\tilde{\phi} (D^l\phi-D^l\tilde{\phi})dx+\int D^l(f_4+f_5) D^l(\phi-\tilde{\phi})dx\\
 &\lesssim \delta_0(\| D^{l +1} (\phi-\tilde{\phi})\|_{L^2}^2+\| D^{l } (\phi-\tilde{\phi})\|_{L^2}^2+\|\phi_x\|_{H^1}^2+\| D^{l }\phi\|_{L^2}^2+\|D^{l } \tilde{\phi}\|_{L^2}^2).   
\end{aligned}
   \end{equation}
Summing $l $ from $0$ to $s$, and noting the smallness of $\delta_0$, one has for $\alpha>0$ arbitrarily small,
   \begin{equation*}
    \frac{1}{2}\frac{d}{dt}\|\phi-\tilde{\phi}\|_{H^s}^2+(\frac{2\bar{v}}{\epsilon}-\alpha)\|\phi-\tilde{\phi}\|_{H^s}^2+\frac{\epsilon}{2}\|(\phi-\tilde{\phi})_x\|_{H^s}^2\lesssim \delta_0\|\tilde{\phi}\|_{H^s}^2,
\end{equation*}
where  the fact that $\|\phi\|_{H^s}\leq\|\phi-\tilde{\phi}\|_{H^s}+\|\tilde{\phi}\|_{H^s}$ is used. Using Gr$\ddot{\text{o}}$nwall's inequality to the above inequality, for the case $\frac{2\bar{v}}{\epsilon}-\alpha<\epsilon$, one obtains
\begin{equation*}
 \| \phi-\tilde{\phi}\|_{H^s}\lesssim \delta_0e^{-(\frac{2\bar{v}}{\epsilon}-\alpha) t}.   
\end{equation*} 
Then for $\bar{\eta} \triangleq \min\{\frac{2\bar{v}}{\epsilon}-\alpha,\epsilon\}$, it holds
\begin{equation*}
  \| \phi-\tilde{\phi}\|_{H^s}\lesssim \delta_0e^{-\bar{\eta} t}.
\end{equation*}
 This completes the proof of Lemma \ref{lem-pd}.
\end{proof}
\vspace{2mm}

Next, we intend to obtain the decay rates of $L^2$-estimates for $(n,w)$ and $(\tilde{n},\tilde{w})$ to nonlinear systems \eqref{nsac-s} and \eqref{ps-s}. To do this, for any $0\leq l\leq s$, introducing the energy
\begin{equation}\label{eslt}
  E_l^s(t)=\|D^ln(t)\|_{H^{s-l}}^2+\|D^l w(t)\|_{H^{s-l}}^2+\|D^l \phi(t)\|_{H^{s-l}}^2,
\end{equation}
and the time-weighted energy functional
\begin{equation}\label{mslt}
  M(t)=\sup_{0\leq\tau\leq t}\left\{(1+\tau)^{\frac{1}{4}}(\|n(\tau)\|_{H^{s}}+\|w(\tau)\|_{H^{s}})\right\}.
\end{equation}
For simplicity, we rewrite the nonlinear systems \eqref{nsac-s}$_{1,2}$ and \eqref{ps-s}$_{1,2}$ into the following systems
\begin{align}
&\partial_tU=\mathcal{L}U+F_x,\quad U|_{t=0}=U_0,\label{nsac-xi-n}\\
&\partial_t\widetilde{U}=\widetilde{\mathcal{L}}\widetilde{U}+\widetilde{F}_x,\quad \widetilde{U}|_{t=0}=U_0,\label{ps-xi-n}
\end{align}
where $U=(n,w)^\top, ~\widetilde{U}=(\tilde{n},\tilde{w})^\top, ~U_0=(n_0,w_0)^\top$ and $F=(0,f_1(n)+f_2(n,w)+f_3(n,\phi))^\top, ~\widetilde{F}=(0,\tilde{f_1}(\tilde{n}))^\top.$
To analysis the large time behavior of the solutions $U$ in frequency space, we introduce the low-high frequency decomposition
\begin{equation}\label{l-h-com}
  U(x,t)=U^{\ell}(x,t)+U^{h}(x,t)=(n^{\ell},w^{\ell})+(n^{h},w^{h}),
\end{equation}
where $U^{\ell}(x,t)=P_1U(x,t), U^{h}(x,t)=P_{\infty}U(x,t)$ and the operator $P_1$ and $P_{\infty}$ are defined by $P_{j}f=\mathscr{F}^{- 1}(\hat\chi_{j}\hat{f}),\quad j=1,\infty$, where $\hat{\chi}_{j}(j=1,\infty)\in C^{\infty}(\mathbb{R}^{3})$, $0\leq \hat\chi_{j}\leq 1$ are cut-off functions defined by
 \begin{equation*}
  \hat\chi_{1}(\xi)=\left\{
  \begin{array}{l}
    1\quad (|\xi|\leq r_{0}),\\
    0\quad (|\xi|\geq R_0),
  \end{array}
  \right.
  \quad\hat\chi_{\infty}(\xi)=1-\hat\chi_{1}(\xi),
 \end{equation*}
 for constants $r_0$ and $R_0$ satisfying $0<r_0<R_0$. For given function $f(x,t)$ and $G^\ell(x,t)=P_1G(x,t)$, it holds
  \begin{equation}\label{pro-gl}
    \|D^k G^\ell*f\|_{L^2}\leq C(1+t)^{-\frac{1}{4}-\frac{k}{2}}\|f\|_{L^1}.
  \end{equation}

\begin{lem}\label{lem-hl}
Under the assumptions of Theorem \ref{thm1}, $(n,w,\phi)$ and $(\tilde{n},\tilde{\phi},\tilde{w})$ are the classical solution of the problem \eqref{nsac-s} and \eqref{ps-s},  then for $0\leq l\leq s$, it holds
  \begin{align}
  &\|D^l(n,w)(t)\|_{H^{s-l}}\lesssim \delta_0(1+t)^{-\frac{1}{4}-\frac{l}{2}},\label{lem-hl0}\\
  &\|D^l(\tilde{n},\tilde{w})(t)\|_{H^{s-l}}\lesssim \delta_0(1+t)^{-\frac{1}{4}-\frac{l}{2}}.\label{lem-ps0}
  \end{align}
\end{lem}
\begin{proof}
First, we prove the estimates \eqref{lem-hl0} for $l=0$.
According to the proof of Proposition \ref{pro-est}, we can obtain
\begin{equation}
  \frac{d}{dt}\|(n,w,\phi)(t)\|_{H^s}^2+\|n_x (t)\|_{H^{s-1}}^2+\|w_x(t)\|_{H^s}^2+\|\phi(t)\|_{H^{s+1}}^2\leq 0.
  \end{equation}
Via the definition of $E_0^s(t)$, one has
  \begin{equation*}
    \frac{d}{dt}E^s_0(t)+\|n_x(t)\|^2_{H^{s-1}}+\|w_x (t)\|_{H^{s-1}}^2+\|\phi(t)\|_{H^{s+1}}^2\leq 0,
  \end{equation*}
which implies that there exists a positive constant $C_3$ such that
  \begin{equation}\label{lem-hs1}
    \frac{d}{dt}E^s_0(t)+C_3E^s_0(t)\lesssim \|(n^{\ell},w^{\ell})(t)\|^2_{L^2},
  \end{equation}
where we have used the fact $\|(n^{h},w^{h})(t)\|_{L^2}\lesssim  \|(n_x,w_x)(t)\|_{L^2}$.
Using Duhamel's principle, for the low frequency part $V^\ell=(n^{\ell},w^{\ell})$, one has
  \begin{equation}\label{lem-hs2}
    V^\ell(x,t)=G^\ell*V_0+\int_0^tG^\ell(t-\tau)*F_x(\tau)d\tau.
  \end{equation}
Using the Plancherel theorem, H$\ddot{\text{o}}$lder's inequality and \eqref{pro-gl}, one has
  \begin{equation}\label{lem-hs3}
    \begin{aligned}
    &\|(n^{\ell},w^{\ell})(t)\|_{L^2}\\
    &\lesssim (1+t)^{-\frac{1}{4}}\|(n_0,w_0)\|_{L^1}+\int_0^t(1+t-\tau)^{-\frac{3}{4}}\|F(\tau)\|_{L^1}d\tau\\
    &\lesssim \delta_0(1+t)^{-\frac{1}{4}}+\int_0^t(1+t-\tau)^{-\frac{3}{4}}(\|n(\tau)\|_{L^2}^2+\|n(\tau)\|_{L^2}\|w_x(\tau)\|_{L^2})d\tau\\
    &\quad+\int_0^t(1+t-\tau)^{-\frac{3}{4}}\|\phi_x(\tau)\|_{L^2}^2d\tau\\
    &\lesssim \delta_0(1+t)^{-\frac{1}{4}}+(M(t))^2\int_0^t(1+t-\tau)^{-\frac{3}{4}}(1+\tau)^{-\frac{1}{2}}d\tau\\
    &\quad+\delta_0^2\int_0^t(1+t-\tau)^{-\frac{3}{4}}e^{-\frac{2\bar{v}\tau}{\epsilon}}d\tau\\
    &\lesssim (1+t)^{-\frac{1}{4}}\left(\delta_0+(M(t))^2\right),
    \end{aligned}
  \end{equation}
where in the second inequality we have used
\begin{equation}
  \|(n,w)(t)\|_{H^s}\leq (1+t)^{-\frac{1}{4}}M(t),\quad \|\phi(t)\|_{H^s}\leq e^{-\frac{\bar{v}t}{\epsilon}}\delta_0.
  \end{equation}
Substituting \eqref{lem-hs3} into \eqref{lem-hs1}, one has
  \begin{equation*}
    \frac{d}{dt}E^s_0(t)+C_3E^s_0(t)\lesssim (1+t)^{-\frac{1}{2}}\left(\delta_0+(M(t))^2\right)^2.
  \end{equation*}
Applying Gr$\ddot{\text{o}}$nwall's inequality to the above inequality, one derives
\begin{equation*}
  \begin{aligned}
      E^s_0(t)&\lesssim e^{-C_3t}E^s_0(0)+ \int_0^te^{-C_3(t-\tau)}(1+\tau)^{-\frac{1}{2}}\left(\delta_0+(M(\tau))^2\right)^2d\tau\\
      &\lesssim (1+t)^{-\frac{1}{2}}\left(\delta_0^2+(M(t))^4\right),
    \end{aligned}
\end{equation*}
which implies that
  \begin{equation}
    (M(t))^2\lesssim \left(\delta_0^2+(M(t))^4\right).
  \end{equation}
Since $\delta_0$ is small enough, it gives rise to
  \begin{equation}\label{lem-hs5}
  M(t)\lesssim \delta_0.
  \end{equation}
It provides
\begin{equation}\label{lem-hs0}
  \|(n,w)(t)\|_{H^s}\lesssim \delta_0(1+t)^{-\frac{1}{4}}.
  \end{equation}

Next, we enhance the decay rate of the estimates for the derivatives.
By \eqref{w2} and \eqref{ph2}, we obtain
\begin{equation}\label{lem-hl2}
  \begin{aligned}
    &\frac{d}{dt}\left(\bar{c}^2\| D^j  n\|_{L^2}^2+\| D^j  w\|_{L^2}^2+\| D^j  \phi\|_{L^2}^2\right)\\    
    &\quad+\bar{\nu}\| D^{j +1}  w \|_{L^2}^2+\frac{2\bar{v}}{\epsilon}\| D^{l } \phi\|_{L^2}^2+\epsilon\| D^{j +1} \phi\|_{L^2}^2\\    
    &\lesssim \delta_0(\| D^j  n\|_{L^2}^2+\| D^{j +1} w\|_{L^2}^2+\| D^{j +1}\phi\|_{L^2}^2+\| D^{j}\phi\|_{L^2}^2).
  \end{aligned}
  \end{equation}
By the same method as in \eqref{n1}, applying the operator $ D^{j-1} P_{\infty}$ to $\eqref{nsac-s}_2$, multiplying the resulting equation by $ D^{j} n^h$, and integrating over $\mathds{R}$, we obtain
\begin{equation}\label{lem-hl3}
\begin{aligned}
  &\frac{\bar{c}^2}{2}\| D^{j}  n^h\|_{L^2}^2-\frac{d}{dt}\int ( D^{j-1}w^h\cdot D^{j}  n^h)dx\\
      &\lesssim C_4(\| D^{j +1} w\|_{L^2}^2+\| D^{j}  w\|_{L^2}^2)+\delta_0(\| D^{j}  n\|_{L^2}^2+\| D^{j +1}\phi\|_{L^2}^2).
\end{aligned}
\end{equation}
Choosing the positive constant $\beta_2$ sufficiently small, then taking $\eqref{lem-hl2}+\beta_2\times(\ref{lem-hl3})$ yields
$$
  \begin{aligned}
  &\frac{d}{dt}\left(\bar{c}^2\| D^j  n\|_{L^2}^2+\| D^j  w\|_{L^2}^2+\| D^{j}  \phi\|_{L^2}^2-\beta_2\int ( D^{j -1}w^h\cdot D^{j}  n^h)dx\right)\\
  &\quad+\bar{\nu}\| D^{j +1}  w \|_{L^2}^2+\frac{2\bar{v}}{\epsilon}\| D^{j} \phi\|_{L^2}^2+\epsilon\| D^{j +1} \phi\|_{L^2}^2+\frac{\bar{c}^2\beta_2}{2}\| D^{j}  n^h\|_{L^2}^2\\    
    &\lesssim (\delta_0+\beta_2C_4)(\| D^{j +1} w\|_{L^2}^2+\| D^{j}  w\|_{L^2}^2)\\    
    &\quad+(\delta_0+\delta_0\beta_2)(\| D^{j}  n\|_{L^2}^2+\| D^{j +1}\phi\|_{L^2}^2+\| D^{j }\phi\|_{L^2}^2).
    \end{aligned}
$$
Since $\beta_2$ and $\delta_0$ are sufficiently small,  one has
\begin{equation}\label{lem-hl4}
  \begin{aligned}
  &\frac{d}{dt}\left(\bar{c}^2\| D^{j}  n\|_{L^2}^2+\| D^{j}  w\|_{L^2}^2+\| D^{j}  \phi\|_{L^2}^2-\beta_2\int ( D^{j -1}w^h\cdot D^{j}  n^h)dx\right)\\
  &\quad+\frac{\bar{\nu}}{2}\| D^{j +1}  w \|_{L^2}^2+\frac{\bar{v}}{\epsilon}\| D^{j} \phi\|_{L^2}^2+\frac{\epsilon}{2}\| D^{j +1} \phi\|_{L^2}^2+\frac{\bar{c}^2\beta_2}{4}\| D^{j}  n^h\|_{L^2}^2\\    
    &\lesssim (\delta_0+\beta_2)(\| D^{j}  w\|_{L^2}^2+\| D^{j}  n\|_{L^2}^2),
    \end{aligned}
\end{equation}
which together with \eqref{a-priori est} yields
\begin{equation}\label{lem-hl5}
\begin{aligned}
&\frac{d}{dt}(\|D^{j}  n\|_{L^2}^2+\|D^{j}  w\|_{L^2}^2+\|D^{j}  \phi\|_{L^2}^2)\\
&\quad+C_5(\|D^{j} n\|_{L^2}^2+\|D^{j +1} w\|_{L^2}^2+\|D^{j} \phi\|_{H^1}^2)\\
&\lesssim (\|D^{j}n^\ell\|_{L^2}^2+\|D^{j}w^\ell\|_{L^2}^2).
\end{aligned}
\end{equation}
Summing up \eqref{lem-hl5} for $j$ from $l$ to $s$ and using the definition of $E_l^s(t)$, one has
\begin{equation}\label{lem-hl-e}
  \frac{d}{dt}E_{l}^s(t)+C_5E_{l}^s(t)\lesssim \| D^{l}(n^{\ell},w^{\ell})\|_{L^2}^2.
  \end{equation}

In what follows, we are going to establish the $L^2$ estimate of $( D^l  n^\ell, D^l  w^\ell)$ with $1\leq l \leq s$. First, similar to the proof of \eqref{lem-hs0}, we also have
\begin{equation}\label{lem-hl6}
  \begin{aligned}
    \|D(n^{\ell},w^{\ell})(t)\|_{L^2}&\lesssim (1+t)^{-\frac{3}{4}}\|(n_0,w_0)\|_{L^1}+\int_0^{\frac t2}(1+t-\tau)^{-\frac{5}{4}}\|F(\tau)\|_{L^1}d\tau\\
    &+\int_{\frac t2}^{t}(1+t-\tau)^{-\frac{3}{4}}\|F_x(\tau)\|_{L^1}d\tau\\
    &\lesssim \delta_0(1+t)^{-\frac{3}{4}},
  \end{aligned}
\end{equation}
which we use the fact that
\begin{align*}
\|F(\tau)\|_{L^1}&\lesssim \|n\|_{L^2}^2+\|n\|_{L^2}\|w_x\|_{L^2}+\|\phi_x\|_{L^2}^2\lesssim \delta_0(1+\tau)^{-\frac{1}{2}},\\
\|F_x(\tau)\|_{L^1}&\lesssim \|n\|_{L^2}\|n_x\|_{L^2}+\|n\|_{L^2}\|w_{xx}\|_{L^2}+\|n_x\|_{L^2}\|w_{x}\|_{L^2}+\|\phi_x\|_{H^1}^2\lesssim \delta_0(1+\tau)^{-\frac{1}{2}}.
\end{align*}
By \eqref{lem-hl-e} and \eqref{lem-hl6}, one gets
\begin{equation*}
  \frac{d}{dt}E_{1}^s(t)+C_5E_{1}^s(t)\lesssim \delta_0(1+t)^{-\frac{3}{2}}.
\end{equation*}
Therefore, it holds
\begin{equation*}
  E_1^s(t)\lesssim \delta_0(1+t)^{-\frac{3}{2}},
  \end{equation*}
which implies
\begin{equation}
  \|D(n,w)(t)\|_{H^{s-1}}\lesssim \delta_0(1+t)^{-\frac{3}{4}}.
  \end{equation}
For second order derivatives, we apply the similar argument.
\begin{equation}\label{lem-hl9}
  \begin{aligned}
    \|D^2(n^{\ell},w^{\ell})(t)\|_{L^2}&\lesssim (1+t)^{-\frac{5}{4}}\|(n_0,w_0)\|_{L^1}+\int_0^{\frac{t}{2}}(1+t-\tau)^{-\frac{7}{4}}\|F(\tau)\|_{L^1}d\tau\\
    &\quad+\int_{\frac{t}{2}}^t(1+t-\tau)^{-\frac{5}{4}}\|F_x(\tau)\|_{L^1}d\tau.
  \end{aligned}
\end{equation}
Note that at this moment, one has
\begin{equation*}
  \|(n,w)(t)\|_{L^2}\lesssim \delta_0(1+t)^{-\frac{1}{4}},\quad \|D(n,w)(t)\|_{H^{s-1}}\lesssim \delta_0(1+t)^{-\frac{3}{4}}.
  \end{equation*}
As a corollary, one obtains
\begin{equation*}
  \begin{aligned}
\|F(\tau)\|_{L^1}&\lesssim \|n\|_{L^2}^2+\|n\|_{L^2}\|w_x\|_{L^2}+\|\phi_x\|_{L^2}^2\lesssim \delta_0(1+\tau)^{-\frac{1}{2}},\\
    \|F_x(\tau)\|_{L^1}&\lesssim \|n\|_{L^2}\|n_x\|_{L^2}+\|n\|_{L^2}\|w_{xx}\|_{L^2}+\|n_x\|_{L^2}\|w_{x}\|_{L^2}+\|\phi_x\|_{H^1}^2\lesssim \delta_0(1+\tau)^{-1}.
  \end{aligned}
  \end{equation*}
Substituting these estimates into \eqref{lem-hl9} yields
\begin{equation*}
  \begin{aligned}
    \|D^2(n^{\ell},w^{\ell})(t)\|_{L^2}&\lesssim \delta_0(1+t)^{-\frac{5}{4}}+\delta_0\int_0^{\frac{t}{2}}(1+t-\tau)^{-\frac{7}{4}}(1+\tau)^{-\frac{1}{2}}d\tau\\
    &\quad+\delta_0\int_{\frac{t}{2}}^t(1+t-\tau)^{-\frac{5}{4}}(1+\tau)^{-1}d\tau\\
    &\lesssim \delta_0(1+t)^{-\frac{5}{4}}.
  \end{aligned}
\end{equation*}
With the help of (\ref{lem-hl-e}), one derives
\begin{equation*}
  \frac{d}{dt}E_{2}^s(t)+C_5E_{2}^s(t)\lesssim \delta_0(1+t)^{-\frac{5}{2}}.
\end{equation*}
Therefore, it holds
\begin{equation*}
  E_2^s(t)\lesssim \delta_0(1+t)^{-\frac{5}{2}},
\end{equation*}
which implies
\begin{equation*}
  \|D^2(n,w)\|_{H^{s-2}}\lesssim \delta_0(1+t)^{-\frac{5}{4}}.
\end{equation*}
The higher order derivatives can be followed in a similar argument and we obtain \eqref{lem-hl0}. The estimates for $D^l(\tilde{n},\tilde{w})$ is clear via the same method as that for $D^l (n,w)$ and we omit the details for simplicity. Thus we finish the proof of Lemma \ref{lem-hl}.
\end{proof}
\vspace{2mm}

Define the energy as follows
\begin{equation}\label{eslt-d}
 \overline{E_k^s}(t)=\|D^k(n-\tilde{n})(t)\|_{H^{s-1-k}}^2+\|D^k (w-\tilde{w})(t)\|_{H^{s-1-k}}^2,
\end{equation}
for any $0\leq k\leq s-1$, we subsequently tend to obtain the decay rates of $L^2$-estimates for $(n-\tilde{n},w-\tilde{w})$. 
\begin{lem}\label{lem-diff0}
Under the assumptions of Theorem \ref{thm1}, $(n,w,\phi)$ and $(\tilde{n},\tilde{\phi},\tilde{w})$ are the classical solution of the problem \eqref{nsac-s} and \eqref{ps-s},  then for $\alpha>0$ arbitrarily small and $0\leq k\leq s-1$,
\begin{align}
  &\|D^k(n-\tilde{n},w-\tilde{w})(t)\|_{H^{s-1-k}}\lesssim \delta_0(1+t)^{-\frac{3}{4}-\frac{k}{2}+\alpha}.\label{lem-diff}
  \end{align}
\end{lem}
\begin{proof}
To begin with, we prove the decay estimates \eqref{lem-diff} for $k=0$. Put
\begin{equation}\label{mslt-d}
  \overline{M}(t)=\sup_{0\leq\tau\leq t}\left\{(1+\tau)^{\frac{3}{4}-\alpha}(\|(n-\tilde{n})(\tau)\|_{H^{s-1}}+\|(w-\tilde{w})(\tau)\|_{H^{s-1}})\right\},
\end{equation}
where $\alpha$ is a small fixed positive constant. It suffices to prove that $\overline{M}_k(t)$ has a uniform time-independent bound. 
Subtracting $\eqref{nsac-s}_{1,2}$ from $\eqref{ps-s}_{1,2}$ yields
\begin{align}\label{nsac-d}
&\left\{\begin{array}{llll}
\displaystyle (n-\tilde{n})_{t}-(w-\tilde{w})_{x}=\frac{\bar{\nu}}{2}(n-\tilde{n})_{xx}-\frac{\bar{\nu}}{2}n_{xx},\\
\displaystyle (w-\tilde{w})_{t}-\bar{c}^2(n-\tilde{n})_{x}=\frac{\bar{\nu}}{2}(w-\tilde{w})_{xx}+\frac{\bar{\nu}}{2}w_{xx}+(f_1 -\tilde{f_1})_x+(f_2+f_3)_x, 
\end{array}\right.
\end{align}
Applying the operator $ D^j~ (j\leq s-1)$ to \eqref{nsac-d}, multiplying the resulting equations by $ \bar{c}^2D^j(n-\tilde{n})$ and $ D^j(w-\tilde{w})$ respectively, and integrating over $\mathds{R}$, one has
\begin{equation*}
\begin{aligned}
 &\frac{1}{2}\frac{d}{dt}\left(\bar{c}^2\| D^j (n\!-\!\tilde{n})\|_{L^2}^2+\| D^j (w\!-\!\tilde{w})\|_{L^2}^2\right)+\frac{\bar{\nu}}{2} (\bar{c}^2\| D^{j+1}(n\!-\!\tilde{n})\|_{L^2}^2+\| D^{j+1}(w\!-\!\tilde{w})\|_{L^2}^2)\\
 &= \frac{\bar{c}^2\bar{\nu}}{2}\int D^jn_xD^{j+1}(n\!-\!\tilde{n})dx-\frac{\bar{\nu}}{2}\int D^jw_xD^{j+1}(w\!-\!\tilde{w})dx-\int D^jf_1 D^{j+1}(w\!-\!\tilde{w})dx\\
 &\quad-\int D^jf_2 D^{j+1}(w\!-\!\tilde{w})dx-\int D^jf_3 D^{j+1}(w\!-\!\tilde{w})dx+\int D^j\tilde{f_1}D^{j+1}(w\!-\!\tilde{w})dx\\
 &\leq \frac{\bar{\nu}}{4} (\bar{c}^2\| D^{j+1}(n\!-\!\tilde{n})\|_{L^2}^2
     +\| D^{j+1}(w\!-\!\tilde{w})\|_{L^2}^2)+C (\|D^jn_x\|_{L^2}^2+\|D^jw_x\|_{L^2}^2)\\
 &\quad +C (\|n\|_{L^\infty}^2\|D^jn\|_{L^2}^2+\|w_x\|_{L^\infty}^2\|D^jn\|_{L^2}^2+\|\tilde{n}\|_{L^\infty}^2\|D^j\tilde{n}\|_{L^2}^2+\|D^j\phi_x\|_{L^2}^2),    
\end{aligned}
   \end{equation*}
which together with \eqref{lem-hl0}-\eqref{lem-ps0} and Lemma \ref{lem-pd} yields
\begin{equation}\label{diff-e1}
\begin{aligned}
 &\frac{d}{dt}\left(\bar{c}^2\| D^j (n\!-\!\tilde{n})\|_{L^2}^2+\| D^j (w\!-\!\tilde{w})\|_{L^2}^2\right)+\frac{\bar{\nu}}{2} (\bar{c}^2\| D^{j+1}(n\!-\!\tilde{n})\|_{L^2}^2+\| D^{j+1}(w\!-\!\tilde{w})\|_{L^2}^2)\\
 &\lesssim  \|D^jn_x\|_{L^2}^2+\|D^jw_x\|_{L^2}^2+ \|n\|_{L^2}\|n_x\|_{L^2}\|D^jn\|_{L^2}^2+\|w_x\|_{H^1}^2\|D^jn\|_{L^2}^2\\
 &\quad +\|\tilde{n}\|_{L^2}\|\tilde{n}_x\|_{L^2}\|D^j\tilde{n}\|_{L^2}^2+\|D^j\phi_x\|_{L^2}^2  \\
 &\lesssim \delta_0(1+t)^{-\frac32-j}.
\end{aligned}
   \end{equation}
Summing $j$ in \eqref{diff-e1} from $0$ to $s-1$, one obtains
\begin{equation*}
\begin{aligned}
 &\frac{d}{dt}\left(\bar{c}^2\|n\!-\!\tilde{n}\|_{H^{s\!-\!1}}^2\!+\!\| w\!-\!\tilde{w}\|_{H^{s\!-\!1}}^2\right)\!+\!\frac{\bar{\nu}}{2} (\bar{c}^2\| (n\!-\!\tilde{n})_x\|_{H^{s\!-\!1}}^2\!+\!\|(w\!-\!\tilde{w})_x\|_{H^{s\!-\!1}}^2)\lesssim \delta_0(1\!+\!t)^{-\frac32},
\end{aligned}
   \end{equation*}
which implies that there exists a positive constant $C_6$ such that
  \begin{equation}\label{diff-e3}
    \frac{d}{dt}\overline{E^s_0}(t)+C_6\overline{E^s_0}(t)\lesssim \|\big((n-\tilde{n})^{\ell},(w-\tilde{w})^{\ell}\big)(t)\|^2_{L^2}+\delta_0(1\!+\!t)^{-\frac32}.
  \end{equation}
Using Duhamel's principle, for the low frequency part $(U-\widetilde{U})^{\ell}=((n-\tilde{n})^{\ell},(w-\tilde{w})^{\ell})$, one has
  \begin{equation}\label{diff-e4}
    (U-\widetilde{U})^\ell(x,t)=(G-\widetilde{G})^\ell*U_0+\int_0^tG^\ell(t-\tau)*F_x(\tau)d\tau-\int_0^t\widetilde{G}^\ell(t-\tau)*\widetilde{F}_x(\tau)d\tau.
  \end{equation}
Using the Plancherel theorem and Lemma \ref{lem-hl}, one gets
\begin{equation}\label{diff-e5}
    \begin{aligned}
    &\|(U-\widetilde{U})^{\ell}(t)\|_{L^2}\\
    &\lesssim (1\!+\!t)^{-\frac{3}{4}}\|U_0\|_{L^1}+\int_0^t(1\!+\!t\!-\!\tau)^{-\frac{3}{4}}\|(f_1 -\tilde{f_1})(\tau)\|_{L^1}d\tau\\
    &\quad +\int_0^t(1\!+\!t\!-\!\tau)^{-\frac{3}{4}}\|(f_2+f_3)(\tau)\|_{L^1}d\tau+\int_0^t(1\!+\!t\!-\!\tau)^{-\frac{3}{4}}\|Df_1 (\tau)\|_{L^1}d\tau\\
    &\lesssim \delta_0(1\!+\!t)^{-\frac{3}{4}}+\delta_0 \overline{M}(t)\int_0^t(1\!+\!t\!-\!\tau)^{-\frac{3}{4}}(1\!+\!\tau)^{-1+\alpha}d\tau\\
    &\quad+\delta_0^2\int_0^t(1\!+\!t\!-\!\tau)^{-\frac{3}{4}}(1\!+\!\tau)^{-1}d\tau+\delta_0^2\int_0^t(1\!+\!t\!-\!\tau)^{-\frac{3}{4}}e^{-\frac{2\bar{v}\tau}{\epsilon}}d\tau\\
    &\lesssim (1+t)^{-\frac{3}{4}+\alpha}\left(\delta_0+\delta_0\overline{M}(t)\right),
    \end{aligned}
  \end{equation}
where in the second inequality, one has used the fact
\begin{equation*}
  \begin{aligned}
&\|(f_1 -\tilde{f_1})(\tau)\|_{L^1}\lesssim (\|n\|_{L^2}+\|\tilde{n}\|_{L^2})\|n-\tilde{n}\|_{L^2}\lesssim \delta_0(1+\tau)^{-1+\alpha}\overline{M}(t),\\
&\|(f_2+Df_1 )(\tau)\|_{L^1}\lesssim (\|n\|_{L^2}\|w_x\|_{L^2}+\|n\|_{L^2}\|n_x\|_{L^2})\lesssim \delta_0^2(1+\tau)^{-1},\\
&\|f_3(\tau)\|_{L^1}\lesssim \|\phi_x\|_{L^2}^2\lesssim \delta_0^2e^{-\frac{2\bar{v}\tau}{\epsilon}}.
  \end{aligned}
  \end{equation*}
Substituting \eqref{diff-e5} into \eqref{diff-e3}, one has
  \begin{equation*}
    \frac{d}{dt}\overline{E^s_0}(t)+C_6\overline{E^s_0}(t)\lesssim (1+t)^{-\frac{3}{2}+2\alpha}\left(\delta_0+\delta_0\overline{M}(t)\right)^2+\delta_0(1+t)^{-\frac{3}{2}}.
  \end{equation*}
Applying Gr$\ddot{\text{o}}$nwall's inequality to the above inequality, one derives 
\begin{equation*}
  \begin{aligned}
      \overline{E^s_0}(t)
      &\lesssim (1+t)^{-\frac{3}{2}+2\alpha}\left(\delta_0^2+\delta_0^2\overline{M}(t)^2\right).
    \end{aligned}
\end{equation*}
Then,  $(\overline{M}(t))^2\lesssim \left(\delta_0^2+\delta_0^2\overline{M}(t)^2\right)$ is obtained, which gives $\overline{M}(t)\lesssim \delta_0 $. This completes the proof of \eqref{lem-diff} for $k=0$.

Now we prove \eqref{lem-diff} for $k\geq 1$. Summing $j$ in \eqref{diff-e1} from $k$ to $s-1$, one has
\begin{equation*}
\begin{aligned}
 &\frac{d}{dt}\left(\bar{c}^2\|D^k(n\!-\!\tilde{n})\|_{H^{s\!-\!1\!-\!k}}^2\!+\!\|D^k(w\!-\!\tilde{w})\|_{H^{s\!-\!1\!-\!k}}^2\right)\!+\!\frac{\bar{\nu}}{2} (\bar{c}^2\|D^k(n\!-\!\tilde{n})_x\|_{H^{s\!-\!1\!-\!k}}^2\!+\!\|D^k(w\!-\!\tilde{w})_x\|_{H^{s\!-\!1\!-\!k}}^2)\\&\lesssim \delta_0(1\!+\!t)^{-\frac32-k},
\end{aligned}
   \end{equation*}
which implies that there exists a positive constant $C_7$ such that
  \begin{equation}\label{diff-e6}
    \frac{d}{dt}\overline{E^s_k}(t)+C_7\overline{E^s_k}(t)\lesssim \|\big(D^k(n-\tilde{n})^{\ell},D^k(w-\tilde{w})^{\ell}\big)(t)\|^2_{L^2}+\delta_0(1\!+\!t)^{-\frac32-k}.
  \end{equation}
Using Lemma \ref{pro-linear} and \eqref{diff-e4}, one gets
\begin{equation*}
    \begin{aligned}
    &\|D\big(W-\widetilde{W}\big)^{\ell}(t)\|_{L^2}\\
    &\lesssim (1\!+\!t)^{-\frac{5}{4}}\|U_0\|_{L^1}+\int_0^t(1\!+\!t\!-\!\tau)^{-\frac{5}{4}}\|(f_1 -\tilde{f_1})(\tau)\|_{L^1}d\tau\\
    &\quad +\int_0^t(1\!+\!t\!-\!\tau)^{-\frac{5}{4}}\|(f_2+f_3)(\tau)\|_{L^1}d\tau+\int_0^{\frac t2}(1\!+\!t\!-\!\tau)^{-\frac{7}{4}}\|f_1 (\tau)\|_{L^1}d\tau\\
    &\quad +\int_{\frac t2}^t(1\!+\!t\!-\!\tau)^{-\frac{5}{4}}\|Df_1 (\tau)\|_{L^1}d\tau.
       \end{aligned}
  \end{equation*}
Note that at this moment, one obtains
\begin{equation*}
  \begin{aligned}
&\|(f_1 -\tilde{f_1})(\tau)\|_{L^1}\lesssim (\|n\|_{L^2}+\|\tilde{n}\|_{L^2})\|n-\tilde{n}\|_{L^2}\lesssim \delta_0(1+\tau)^{-1+\alpha},
  \end{aligned}
\end{equation*}
which together with Lemma \ref{lem-hl} implies
\begin{equation*}
    \begin{aligned}
    &\|D\big(W-\widetilde{W}\big)^{\ell}(t)\|_{L^2}\lesssim \delta_0(1+t)^{-\frac54+\alpha}.       
\end{aligned}
\end{equation*}
Substituting this into \eqref{diff-e6} with $k=1$ yields
\begin{equation*}
    \frac{d}{dt}\overline{E^s_1}(t)+C_7\overline{E^s_1}(t)\lesssim \delta_0(1\!+\!t)^{-\frac52+2\alpha}.
  \end{equation*}
Then it holds
\begin{equation*}
    \overline{E^s_1}(t)\lesssim \delta_0(1\!+\!t)^{-\frac52+2\alpha},
  \end{equation*}
which implies that
\begin{equation*}
  \|D(n-\tilde{n},w-\tilde{w})(t)\|_{H^{s-2}}\lesssim \delta_0(1+t)^{-\frac{5}{4}+\alpha}.
  \end{equation*}
The higher order derivatives can be followed in a similar argument and we obtain \eqref{lem-diff} and  complete the proof of Lemma \ref{lem-diff0}.
\end{proof}

Finally, we tend to obtain the decay rates of $L^1$-estimates for $(n-\tilde{n},w-\tilde{w},\phi-\tilde{\phi})$. 
\begin{lem}\label{lem-diff1}
Under the assumptions of Theorem \ref{thm2}, $(n,w,\phi)$ and $(\tilde{n},\tilde{\phi},\tilde{w})$ are the classical solution of the problem \eqref{nsac-s} and \eqref{ps-s}, then  for $\alpha>0$ arbitrarily small and $0\leq k\leq s-1$,
\begin{align}
  &\|D^k(n-\tilde{n},w-\tilde{w})(t)\|_{W^{s-1-k,1}}\lesssim \delta_0(1+t)^{-\frac{1}{2}-\frac{k}{2}+\alpha}+e^{-\frac{\bar{c}^2}{\bar{\nu}}t}\|U_0\|_{W^{s-1,1}},\label{diff1-1}\\
  &\|(\phi-\tilde{\phi})(t)\|_{W^{s,1}}\lesssim (\delta_0^2+\|\phi_0\|_{L^1})e^{-\bar{\eta} t}.\label{diff1-2}
  \end{align}
where $\bar{\eta}\triangleq \min\{\frac{2\bar{v}}{\epsilon}-\alpha,\epsilon\}.$
\end{lem}
\begin{proof}
 From \eqref{nsac-xi-n} and \eqref{ps-xi-n}, one has 
\begin{align}
\displaystyle & D^jU(x,t)=D^jG*U_0+\int_0^tD^j(G(t-\tau)*F_x(\tau))d\tau,\\
\displaystyle & D^j\widetilde{U}(x,t)=D^j\widetilde{G}*U_0+\int_0^tD^j(\widetilde{G}(t-\tau)*\widetilde{F}_x(\tau))d\tau.
\end{align}
Then it holds
\begin{equation}
\begin{aligned}
 D^j(U-\widetilde{U})(x,t)=&D^j(G-\widetilde{G})*U_0+\int_0^tD^j((G-\widetilde{G})(t-\tau)*Df_1 (\tau))d\tau\\
 &+\int_0^tD^j(G(t-\tau)*D(f_2+f_3)(\tau))d\tau\\
 &+\int_0^tD^j(\widetilde{G}(t-\tau)*D(f_1 -\tilde{f_1})(\tau))d\tau.
\end{aligned}  
\end{equation}
By Lemma \ref{pro-l1}, for $j\leq s-1$,
\begin{eqnarray}\label{diff1-3}
&&\|D^j(U-\widetilde{U})(t)\|_{L^1}\notag\\
&&\lesssim (1\!+\!t)^{-\frac{1+j}{2}}\|U_0\|_{L^1}+e^{-\frac{\bar{c}^2}{\bar{\nu}}t}\|U_0\|_{W^{j,1}}\notag\\
&&\quad+\int_0^{\frac t2} \left((1\!+\!t\!-\!\tau)^{-\frac{2+j}{2}}\|f_1 \|_{L^1}+e^{-\frac{\bar{c}^2(t\!-\tau\!)}{\bar{\nu}}}\|f_1 \|_{W^{j+1,1}}\right)d\tau\notag\\
&&\quad+\int_{\frac t2}^t \left((1\!+\!t\!-\!\tau)^{-\frac{1}{2}}\|D^{j+1}f_1 \|_{L^1}+e^{-\frac{\bar{c}^2(t\!-\!\tau)}{\bar{\nu}}}\|D^jf_1 \|_{W^{1,1}}\right)d\tau\notag\\
&&\quad+\int_0^{\frac t2} \left((1\!+\!t\!-\!\tau)^{-\frac{1+j}{2}}(\|f_2\|_{L^1}\!+\!\|f_3\|_{L^1})+e^{-\frac{\bar{c}^2(t\!-\!\tau)}{\bar{\nu}}}(\|f_2\|_{W^{j+1,1}}\!+\!\|f_3\|_{W^{j+1,1}})\right)d\tau\notag\\
&&\quad+\int_{\frac t2}^t \left((1\!+\!t\!-\!\tau)^{-\frac{1}{2}}(\|D^{j}f_2\|_{L^1}\!+\!\|D^{j}f_3\|_{L^1})+e^{-\frac{\bar{c}^2(t\!-\!\tau)}{\bar{\nu}}}(\|D^jf_2\|_{W^{1,1}}\!+\!\|D^jf_3\|_{W^{1,1}})\right)d\tau\notag\\
&&\quad+\int_0^{\frac t2} \left((1\!+\!t\!-\!\tau)^{-\frac{1+j}{2}}\|f_1 \!-\!\tilde{f_1}\|_{L^1}+e^{-\frac{\bar{c}^2(t\!-\!\tau)}{\bar{\nu}}}(\|f_1 \|_{W^{j+1,1}}\!+\!\|\tilde{f_1}\|_{W^{j+1,1}})\right)d\tau\notag\\
&&\quad+\int_{\frac t2}^t \left((1\!+\!t\!-\!\tau)^{-\frac{1}{2}}\|D^{j}(f_1 \!-\!\tilde{f_1})\|_{L^1}+e^{-\frac{\bar{c}^2(t\!-\!\tau)}{\bar{\nu}}}(\|D^jf_1 \|_{W^{1,1}}\!+\!\|D^j\tilde{f_1}\|_{W^{1,1}})\right)d\tau.
\end{eqnarray}  
By Lemmas \ref{lem-pd}-\ref{lem-diff0}, it is easy to verify that
\begin{equation*}
\begin{aligned}
& \|D^jf_1 (t)\|_{L^{1}}\lesssim \sum_{k=0}^j\|D^{j-k}nD^kn\|_{L^{1}}\lesssim \sum_{k=0}^j\|D^{j-k}n\|_{L^2}\|D^kn\|_{L^{2}}\lesssim \delta_0^2 (1\!+\!t)^{-\frac{1+j}{2}},\\
& \|D^jf_2 (t)\|_{L^{1}}\lesssim \sum_{k=0}^j\|D^{j-k}nD^kw_x\|_{L^{1}}\lesssim \delta_0^2 (1\!+\!t)^{-\frac{2+j}{2}},\\
& \|D^jf_3 (t)\|_{L^{1}}\lesssim \sum_{k=0}^j\|D^{j-k}(\frac{\phi}{n+\bar{v}})D^k\phi_x\|_{L^{1}}\lesssim \delta_0^2 e^{-\eta t},\\
& \|D^j(f_1\!-\!\tilde{f_1}) (t)\|_{L^{1}}\lesssim \sum_{k=0}^j\|D^{j-k}(n+\tilde{n})D^k(n-\tilde{n})\|_{L^{1}}\lesssim \delta_0^2 (1\!+\!t)^{-\frac{2+j}{2}+\alpha}.
 \end{aligned}  
\end{equation*}
Combining these estimates into \eqref{diff1-3} yields 
\begin{equation}\label{diff1-4}
\begin{aligned}
 \|D^j(U-\widetilde{U})\|_{L^1}(t)
 &\lesssim (1\!+\!t)^{-\frac{1+j}{2}+\alpha}\delta_0+e^{-\frac{\bar{c}^2}{\bar{\nu}}t}\|U_0\|_{W^{j,1}}.
\end{aligned}  
\end{equation}
Summing up \eqref{diff1-4} for $j$ from $0$ to $s-1$ gives \eqref{diff1-1}.

Finally, we tend to prove \eqref{diff1-2}. From \eqref{nsac-s}, \eqref{ps-s} and \eqref{phi-ps}, we have
\begin{equation}\label{diff1-5}
\begin{aligned}
 D^l(\phi-\tilde{\phi})&=e^{-\frac{2\bar{v}}{\epsilon}t}D^l(H(x,t)*\phi_0(x))-e^{-\epsilon t}D^l(H(x,t)*\phi_0(x))\\
 &\quad+\int_0^te^{-\frac{2\bar{v}}{\epsilon}(t-\tau)}D^l(H(x,t-\tau)*(f_4+f_5)(\tau))d\tau,
\end{aligned}  
\end{equation}
where $H(x,t)=\frac{1}{\sqrt{4\pi\epsilon t}}e^{-\frac{x^2}{4\epsilon t}}$.
By Lemmas \ref{pro-l1}-\ref{lem-pd}, we obtain
\begin{equation}\label{diff1-6}
\begin{aligned}
& \|D^l(\phi-\tilde{\phi})(t)\|_{L^1}\\
&\lesssim e^{-\eta t}\|\phi_0\|_{L^1}+\int_0^te^{-\frac{2\bar{v}(t\!-\tau\!)}{\epsilon}}(1\!+\!t\!-\!\tau)^{-\frac{l}{2}}\left(\|f_4\|_{L^{1}}+\|f_5\|_{L^{1}}\right)d\tau\\
&\lesssim e^{-\eta t}\|\phi_0\|_{L^1}+\int_0^te^{-\frac{2\bar{v}(t\!-\tau\!)}{\epsilon}}(1\!+\!t\!-\!\tau)^{-\frac{l}{2}}\left(\|n\|_{H^{1}}\|\phi\|_{H^{1}}+\|\phi\|_{H^{1}}^2\right)d\tau\\
&\lesssim e^{-\eta t}\|\phi_0\|_{L^1}+\delta_0^2\int_0^te^{-\frac{2\bar{v}(t\!-\tau\!)}{\epsilon}}(1\!+\!t\!-\!\tau)^{-\frac{l}{2}}(1\!+\!\tau)^{-\frac{1}{4}}e^{-\frac{2\bar{v}\tau}{\epsilon}+\alpha\tau}d\tau\\
&\quad+\delta_0^2\int_0^te^{-\frac{2\bar{v}(t\!-\tau\!)}{\epsilon}}(1\!+\!t\!-\!\tau)^{-\frac{l}{2}}e^{-\frac{4\bar{v}\tau}{\epsilon}+2\alpha\tau}d\tau\\
&\lesssim e^{-\eta t}\|\phi_0\|_{L^1}+\delta_0^2e^{-\frac{2\bar{v}t}{\epsilon}+\alpha t}\int_0^te^{-\alpha(t\!-\tau\!)}(1\!+\!t\!-\!\tau)^{-\frac{l}{2}}(1\!+\!\tau)^{-\frac{1}{4}}d\tau\\
&\quad+\delta_0^2e^{-\frac{2\bar{v}t}{\epsilon}+\alpha t}\int_0^te^{-\frac{2\bar{v}\tau}{\epsilon}+\alpha \tau}e^{-\alpha(t\!-\tau\!)}(1\!+\!t\!-\!\tau)^{-\frac{l}{2}}d\tau\\
&\lesssim e^{-\eta t}\|\phi_0\|_{L^1}+\delta_0^2e^{-\frac{2\bar{v}t}{\epsilon}+\alpha t}\\
&\lesssim (\delta_0^2+\|\phi_0\|_{L^1})e^{-\bar{\eta} t},
\end{aligned}  
\end{equation}
where $\eta=\min\{-\frac{2\bar{v}t}{\epsilon},\epsilon\}$, $\bar{\eta}=\min\{-\frac{2\bar{v}t}{\epsilon}+\alpha,\epsilon\}$, and $\alpha$ is a small positive constant. Sum up \eqref{diff1-6} for $l$ from $0$ to $s$ yields the inequality \eqref{diff1-2}.
This completes the proof.
\end{proof}

\subsection{The proof of Theorems \ref{thm1}-\ref{thm2}}
\hspace{2em}With the a priori estimates in Section \ref{se2} and decay estimates in Section \ref{Sec4} at hand, we are now in a position to complete the proof of the theorems in our paper.
\vspace{2mm}

\noindent{\it\textbf{Proof of Theorem \ref{thm1}.}\ }
  By Proposition \ref{pro-est}, it gives
    \begin{equation}
\|(n,w,\phi)(t)\|_{H^s}^2+\int_0^t\left(\|n_x\|_{H^{s-1}}^2+\|w_x\|_{H^{s}}^2+\|\phi\|_{H^{s+1}}^2\right)d\tau\lesssim \delta_0^2.
  \end{equation}
By choosing the initial data $\delta_0$ sufficiently small such that $C\delta_0^2 \leq \frac{1}{4}\delta^2$, we could close the a-priori assumption \eqref{a-priori est}. Then, based on the continuous argument, the global existence of solution $(v,u,\varphi)$ and the estimates \eqref{con2} is obtained.
Combining \eqref{pertub} with Lemma \ref{lem-pd} and Lemma \ref{lem-hl}, it is easy to yields \eqref{lt2}. This completes the proof of Theorem \ref{thm1}.\endproof
\vspace{2mm}

\noindent{\it\textbf{Proof of Theorem \ref{thm2}.}\ }
Due to \eqref{pertub} and \eqref{pertub2}, one has
  \begin{equation}\label{thm2-p1}
   v-\tilde{v}=n-\tilde{n}, \quad u-\tilde{u}=w-\tilde{w}, \quad \varphi-\tilde{\varphi}=\phi-\tilde{\phi},
  \end{equation}
which together with Lemma \ref{lem-pd} and Lemma \ref{lem-diff0} yields \eqref{lt3}. Lemma \eqref{lem-diff1} and \eqref{thm2-p1} yields \eqref{lt4} and completes the proof of Theorem \ref{thm2}.
\endproof

\section*{Acknowledgements}
\hspace{2em}Qiaolin He and Yazhou Chen acknowledge support from  National Natural Sciences Foundation of China  (NSFC) (No. 12371434). Qiaolin He acknowledges support from the National key R \& D Program of China (No.2022YFE03040002).
· Xiaoding Shi acknowledges support from  National Natural Sciences Foundation of China  (NSFC) (No. 12171024).

\end{document}